\documentclass{amsart}

\usepackage{amssymb}
\usepackage{eepic}
\usepackage{color}
\usepackage[linktocpage=true]{hyperref}

\allowdisplaybreaks

\numberwithin{equation}{section}

\newtheorem{theorem}{Theorem}[section]
\newtheorem{lemma}[theorem]{Lemma}

\newtheorem{corollary}[theorem]{Corollary}
\newtheorem{remark}[theorem]{Remark}

\newtheorem{TheoA}{Theorem A}

\newtheorem{TheoB}{Theorem B}

\newcommand{\C}{\mathbb{C}}

\newcommand{\summ}{\sum\nolimits}

\def\1{\mathbf{1}}
\def\H{\mathcal{H}}
\def\E{\mathcal{E}}

\def\M{\mathcal{M}}

\newcommand{\dem}{\noindent {\bf Proof. }}

\newcommand{\fin}{\hspace*{\fill} $\square$ \vskip0.2cm}

\addtocounter{tocdepth}{0}

\begin{document}

\null

\vskip-40pt

\null

\title[Asymmetric Doob inequalities]{Algebraic Davis decomposition \\ and asymmetric Doob inequalities}

\author[Hong, Junge, Parcet]
{Guixiang Hong, Marius Junge, Javier Parcet}

\maketitle

\null

\vskip-45pt

\null

\begin{abstract}
In this paper we investigate asymmetric forms of Doob maximal inequality. The asymmetry is imposed by noncommutativity. Let $(\M,\tau)$ be a noncommutative probability space equipped with a weak-$*$ dense filtration of von Neumann subalgebras $(\M_n)_{n \ge 1}$. Let $\E_n$ denote the corresponding family of conditional expectations. As an illustration for an asymmetric result, we prove that for $1 < p < 2$ and $x \in L_p(\M,\tau)$ one can find $a, b \in L_p(\M,\tau)$ and contractions $u_n, v_n \in \M$ such that $$\E_n(x) = a u_n + v_n b \quad \mbox{and} \quad \max \big\{ \|a\|_p, \|b\|_p \big\} \le c_p \|x\|_p.$$ Moreover, it turns out that $a u_n$ and $v_n b$ converge in the row/column Hardy spaces $\H_p^r(\M)$ and $\H_p^c(\M)$ respectively. In particular, this solves a problem posed by Defant and Junge in 2004. In the case $p=1$, our results establish a noncommutative form of Davis celebrated theorem on the relation between martingale maximal and square functions in $L_1$, whose noncommutative form has remained open for quite some time. Given $1 \le p \le 2$, we also provide new weak type maximal estimates, which imply in turn left/right almost uniform convergence of $\E_n(x)$ in row/column Hardy spaces. This improves the bilateral convergence known so far. Our approach is based on new forms of Davis martingale decomposition which are of independent interest, and an algebraic atomic description for the involved Hardy spaces. The latter results are new even for commutative von Neumann algebras.
\end{abstract}

\addtolength{\parskip}{+1ex}

\renewcommand{\theequation}{AD$_{L_p}$}
\addtocounter{equation}{-1}

\null

\vskip-25pt

\null

\section*{{\bf Introduction}}

Doob maximal inequality is a corner stone in harmonic analysis, probability and ergodic theory. Its noncommutative form is central in noncommutative harmonic analysis and quantum probability. Cuculescu established in \cite{Cu} the noncommutative endpoint estimate for $p=1$ of Doob's inequality. Given $(\M,\tau)$ a noncommutative probability space, let $\E_n$ denote the conditional expectantions associated to a given weak-$*$ dense filtration $(\M_n)_{n \ge 1}$. Given $x \in L_1(\M)_+$ and $\lambda>0$, Cuculescu constructed projections $q_\lambda \in \M$ satisfying $$q_\lambda \E_n(x) q_\lambda \le \lambda \qquad \mbox{and} \qquad \tau \big( \1 - q_\lambda \big) \le \frac{1}{\lambda} \|x\|_1.$$ Unfortunately, Marcinkiewicz interpolation with the other (obvious) endpoint is by no means trivial. Due to the lack of pointwise suprema after quantization, it first required to understand how noncommutative $L_p$ norms of maximal functions should be described. This was achieved by Pisier using sophisticated tools from operator space theory \cite{Pis98}. Then, the expected interpolation result was proved by Junge/Xu in 2007 for positive cones \cite{JuXu07}. A few years earlier, the second-named author had found a direct more elaborated argument from Hilbert module theory \cite{Jun02}. Given $p>1$ and $x \in L_p(\M)$, the noncommutative form of Doob maximal $L_p$ inequality provides operators $a, b \in L_{2p}(\M)$ and $w_n \in \M$ satisfying $$\E_n(x) = a w_n b \qquad \mbox{and} \qquad \|a\|_{2p} \Big( \sup_{n \ge 1} \|w_n\|_\M \Big) \|b\|_{2p} \le c_p \|x\|_p.$$ The results above reduce to Doob's original formulation for commutative algebras. 

As we shall see, the spaces above have a symmetric nature. Noncommutativity allows however to conjecture natural asymmetric forms of these inequalities, which all collapse into one inequality for abelian algebras. The row/column-valued $L_p$ spaces ---the most asymmetric ones--- are omnipresent in operator space theory and quantum probability. Just to mention some examples, noncommutative Khintchine or Burkholder-Gundy inequalities \hskip-0.7pt \cite{Lu,PR,PX1}, as well as several noncommutative forms of Littlewood-Paley theory \cite{JLMX,JMP} \, precise row/column spaces. Certain free variants of these inequalities have also a great impact in Grothendieck's theorem for operator spaces \cite{PS,Xu}. Other more subtle asymmetries were studied in \cite{JuPa10} with applications in operator space $L_p$ embedding theory \cite{JuPa08}. In the particular context of noncommutative maximal inequalities, almost everywhere convergence is replaced by almost uniform convergence and row/column asymmetric estimates yield left/right a.u. convergence \cite[Proposition 5.1]{DeJu04}, less restrictive than what the symmetric ones provide. The row/column $L_p(\ell_\infty)$ spaces have also played a role in noncommutative BMO theory. All of this motivates further research.

Motivated by a question of Gilles Pisier, the second-named author found the first asymmetric forms of Doob inequality in his paper \cite{Jun02}. Namely, given $0 \le \theta \le 1$ the following estimate holds for $x \in L_p(\M)$ 
\begin{equation} \label{Asymmetric1}
\big\| (\E_n(x))_{n \ge 1} \big\|_{L_p(\M; \ell_\infty^\theta)} \le c_{p,\theta} \|x\|_{L_p(\M)} \quad \mbox{when} \quad p > 2 \max \big\{ \theta, 1 - \theta \big\}.
\end{equation}
$L_p(\M;\ell_\infty^\theta)$ is the subspace of sequences in $L_p(\M)$ with quasi-norm $$\big\| (x_n)_{n \ge 1} \big\|_{L_p(\M;\ell_\infty^\theta)} = \inf \Big\{ \|a\|_{\frac{p}{1-\theta}} \Big( \sup_{n \ge 1} \|w_n\|_\infty \Big) \|b\|_{\frac{p}{\theta}} \, \big| \ x_n = a w_n b \mbox{ for } n \ge 1\Big\}.$$ The infimum 
\renewcommand{\theequation}{AD$_{\H_p}$}
\addtocounter{equation}{-1}
runs over all possible factorizations of $(x_n)_{n \ge 1}$ in the form $x_n = a w_n b$ with $(a,b) \in L_{p/(1-\theta)}(\M) \times L_{p/\theta}(\M)$ and $(w_n)_{n \ge 1}$ uniformly bounded in $\M$. The symmetric Doob inequality corresponds to $\theta=1/2$ and the corresponding space is denoted $L_p(\M; \ell_\infty)$. Other significant cases are given  by the row/column spaces $L_p(\M; \ell_\infty^r)$ and $L_p(\M; \ell_\infty^c)$ which correspond to $\theta = 0,1$ respectively. Let us note that the triangle inequality may fail unless $1 - p/2 \le \theta \le p/2$ ---equivalently $p \ge 2 \max \{\theta, 1 - \theta\}$--- and these spaces form a natural interpolation scale in this range \cite{DeJu04,JuPa10}. Although \eqref{Asymmetric1} is fully satisfactory for $p>2$, a counterexample in \cite{DeJu04} disproved the asymmetric inequality for $p < 2 \max \{ \theta, 1 - \theta \}$. This implies that row/column estimates fail for $p < 2$. Fortunately, this is not the end of the story and other asymmetric Doob inequalities might hold. Indeed, given $1 \le p \le 2$ and recalling $\H_p^r(\M) + \H_p^c(\M) \subset L_p(\M)$ from the noncommutative Burkholder-Gundy inequalities \cite{PX1}, the best we could hope for is   
\begin{equation} \label{Asymmetric2}
\begin{array}{rcl} \big\| (\E_n (x))_{n \ge 1} \big\|_{L_p(\M; \ell_\infty^r)} \!\!\!\! & \le & \!\!\!\! c_p \|x\|_{\H_p^r(\M)}, \\
\big\| (\E_n (x))_{n \ge 1} \big\|_{L_p(\M; \ell_\infty^c)} \!\!\!\! & \le & \!\!\!\! c_p \|x\|_{\H_p^c(\M)}. \end{array}
\end{equation}
The row/column Hardy spaces $\H_p^r(\M)$ and $\H_p^c(\M)$ are the completion of finite $L_p$ martingales with respect to the $p$-norm of their (row/column) martingale square functions. This suggests a control of row/column maximal functions by row/column square functions in the spirit of Davis fundamental theorem \cite{Da}. In the symmetric case of $\H_p(\M)$, it holds for $p > 1$ \cite{Jun02,PX1} and fails for $p=1$ \cite{JX0}. Unfortunately it seems that \eqref{Asymmetric2} is too good to be true ---see below--- but we may find their closest substitutes. Our first result establishes weak type forms of \eqref{Asymmetric2} and also strong forms after arbitrary small perturbations of the asymmetries.  
 
There are two ways to provide weak forms of the space $L_p(\M;\ell_\infty)$. One is as an amalgamated space $L_{p,\infty}(\M; \ell_\infty) = L_{2p,\infty}(\M) \ell_\infty(\M) L_{2p,\infty}(\M)$ in the spirit of \cite{JuPa10}, which allows asymmetric generalizations in an obvious way. Alternatively the (weaker) space $\Lambda_{p,\infty}(\M,\ell_\infty)$ is defined as the sequences $(x_n)_{n \ge 1}$ in $L_{p,\infty}(\M)$ satisfying that $$\big\| (x_n)_{n \ge 1} \big\|_{\Lambda_{p, \infty}(\M;\ell_\infty)} = \sup_{\lambda > 0} \inf_{q \in \M_\pi} \Big\{ \lambda \tau \big( (\1 - q) \big)^{\frac{1}{p}} \, \big| \ \|q x_n q\|_\infty \le \lambda \mbox{ for all } n \ge 1 \Big\}$$ is finite. Here $\M_\pi$ stands for the projection lattice in $\M$. This definition is inspired by Cuculescu's construction and $L_{p,\infty}(\M; \ell_\infty) \subset \Lambda_{p,\infty}(\M; \ell_\infty)$. The column space is  determined by $$\big\| (x_n)_{n \ge 1} \big\|_{\Lambda_{p, \infty}(\M;\ell_\infty^c)} \hskip1pt = \hskip1pt \sup_{\lambda > 0} \inf_{q \in \M_\pi} \Big\{ \lambda \tau \big( (\1 - q) \big)^{\frac{1}{p}} \, \big| \ \| x_n q\|_\infty \le \lambda \mbox{ for all } n \ge 1 \Big\}.$$ This is finite iff $(x_n^* x_n)_{n \ge 1} \in \Lambda_{p/2,\infty}(\M; \ell_\infty)$. Take adjoints to define the row space. 

\begin{TheoA}
Let $(\M,\tau)$ be a noncommutative probability space and let $\E_n$ denote the conditional expectations associated to a weak-$*$ dense filtration $(\M_n)_{n \ge 1}$ of von Neumann sulbalgebras. Then, the following inequalities hold:
\begin{itemize}
\item[i)] Given $1 \le p \le 2$ and $x \in \H_p^c(\M)$ $$\big\| (\E_n(x))_{n \ge 1} \big\|_{L_p(\M;\ell_\infty^\theta)} \le c_{p,\theta} \|x\|_{\H_p^c(\M)}$$ provided $1 - p/2 < \theta < 1$. The same holds for $x \in \H_p^r(\M)$ and $0 < \theta < p/2$.

\vskip3pt 

\item[ii)] Given $1 \le p \le 2$ and $x \in \H_p^c(\M)$ $$\big\| (\E_n(x))_{n \ge 1} \big\|_{\Lambda_{p,\infty}(\M;\ell_\infty^c)} \le c_p \|x\|_{\H_p^c(\M)}.$$
Similarly, the row analog $(\E_n)_{n \ge 1}: \H_p^r(\M) \to \Lambda_{p,\infty}(\M;\ell_\infty^r)$ is also bounded.
\end{itemize}
\end{TheoA}

Theorem A gets very close to inequalities \eqref{Asymmetric2} ---see Theorem B for related inequalities--- and according to \cite{DeJu04,JX0} we conjecture that Theorem A is best possible in our restrictions on $0 \le \theta \le 1$, see Remark \ref{Optimalidad}. Actually, Theorem A solves the mystery around the noncommutative Davis theorem and yields the following result for any $0 < \theta \neq \frac12 < 1$ $$(\E_n)_{n \ge 1}: \H_1(\M) \to L_1(\M; \ell_\infty^\theta) + L_1(\M;\ell_\infty^{1-\theta}).$$ The symmetric case $\theta = \frac12$ was disproved in \cite{JX0} but Theorem A shows it works for arbitrary small asymmetries. \hskip-2pt Theorem Aii also provides weak estimates for $\theta = 0,1$.

A crucial difficulty in the proof is that we may not take direct advantage of the positivity-preserving nature of conditional expectations, as it happens in previous results \cite{DeJu04,Jun02,JuXu07}. Our proof rests on two crucial points. We first decompose the column Hardy space $\H_p^c(\M) = h_p^c(\M) + h_p^d(\M)$ into the column conditioned Hardy space $h_p^c(\M)$ and the diagonal Hardy space $h_p^d(\M)$, precise definitions will be given below in the body of the paper. This result is known as the noncommutative Davis decomposition, independently discovered by Junge/Mei and Perrin \cite{JuMe10,Per0} and subsequently improved in \cite{JuPe14,Per} with a better diagonal term $h_p^{1_c}(\M)$. The second ingredient is an instrumental \lq algebraic atomic\rq${}$ description of these spaces from \cite{Per}. The combination of these two results produces a description of $\H_p^c(\M)$ which we call \emph{algebraic Davis decomposition} in this paper.   

A sequence $(x_n)_{n \ge 1}$ of $\tau$-measurable operators converges to $0$ $\tau$-almost uniformly when there is a sequence of projections $(p_k)_{k \ge 1}$ in $\M$ satisfying $\lim_k \tau( \1 - p_k) = 0$ and $\lim_n \|x_n p_k\|_\infty = 0$ for all $k \ge 1$. Theorem Aii implies that $(\E_n(x))_{n \ge 1}$ converges $\tau$-a.u. to $x$ for every $x \in \H_p^c(\M)$ and any $1 \le p \le 2$. Row analogs also apply and refine the (weaker) $\tau$-a.u. bilateral convergence results in \cite{Jun02}. 

\renewcommand{\theequation}{DD$_{pw}$}
\addtocounter{equation}{-1}

Stronger asymmetric Doob maximal estimates follow by stretching our approach to produce finer algebraic Davis type decompositions. More precisely, according to Theorem A and the noncommutative  Burkholder-Gundy inequalities, it is clear that every $x$ in $L_p(\M)$ can be written as $x = x_r + x_c$ with $$\max \Big\{ \big\| (\E_n(x_r))_{n \ge 1} \big\|_{\Lambda_{p,\infty}(\M;\ell_\infty^r)}, \big\| (\E_n(x_c))_{n \ge 1} \big\|_{\Lambda_{p,\infty}(\M;\ell_\infty^c)} \Big\} \, \le \, c_p \|x\|_{L_p(\M)}.$$ Similar decompositions apply for the strong  inequalities in Theorem Ai. Can we find a better decomposition $x = x_r + x_c$ to prove the inequality above for row/column $L_p(\ell_\infty)$ spaces instead of their weak analogs? In this paper we will introduce new families of spaces 
\begin{equation} \label{AlgDavis}
\underbrace{h_{pw}^r(\M) + h_{pw}^{1_r}(\M)}_{\H_{pw}^r(\M)} \quad \mbox{and} \quad \underbrace{h_{pw}^c(\M) + h_{pw}^{1_c}(\M)}_{\H_{pw}^c(\M)}
\end{equation}
for $w \ge 2$, so that the spaces corresponding to the parameter $w=2$ recover the row/column Hardy spaces considered so far. The key to solve the question above is a new algebraic Davis decomposition which refines the ones in \cite{JuMe10,JuPe14,Per0, Per}. We think it is of independent interest. In the following result we include this Davis decomposition and the strong type inequality which answers our question. 

\numberwithin{equation}{section}

\begin{TheoB}
Let $(\M,\tau)$ be a noncommutative probability space and let $\E_n$ denote the conditional expectations associated to a weak-$*$ dense filtration $(\M_n)_{n \ge 1}$ of von Neumann sulbalgebras. Then, the following results hold:
\begin{itemize}
\item[i)] Given $1 < p < 2$ with $1/p = 1/w + 1/s$, we find $$L_p(\M) \simeq \H_{pw}^r(\M) + \H_{pw}^c(\M) \quad \mbox{provided} \quad w, s \ge 2.$$ Moreover, we have continuous inclusions $\H_{pw}^\dag(\M) \subset \H_p^\dag(\M)$ for $\dag = r, c$.

\vskip3pt 

\item[ii)] Given $1 < p < 2$, the inequalities below hold for any $w > 2$ 
\begin{eqnarray*}
\big\| (\E_n(x))_{n \ge 1} \big\|_{L_{p}(\M;\ell_\infty^r)} \!\!\!\! & \le & \!\!\!\! c_{p,w} \|x\|_{\H_{pw}^r(\M)}, \\
\big\| (\E_n(x))_{n \ge 1} \big\|_{L_{p}(\M;\ell_\infty^c)} \!\!\!\! & \le & \!\!\!\! c_{p,w} \|x\|_{\H_{pw}^c(\M)}.
\end{eqnarray*}
\noindent In particular, given $x \in L_p(\M)$ may write $x = x_r + x_c$ with $$\hskip30pt \max \Big\{ \big\| (\E_n(x_r))_{n \ge 1} \big\|_{L_{p}(\M;\ell_\infty^r)}, \big\| (\E_n(x_c))_{n \ge 1} \big\|_{L_{p}(\M;\ell_\infty^c)} \Big\} \, \le \, c_p \|x\|_{L_p(\M)}.$$ Moreover, we have $x_r \in \H_{pw}^r(\M) \subset \H_p^r(\M)$ and $x_c \in \H_{pw}^c(\M) \subset \H_p^c(\M)$.
\end{itemize}
\end{TheoB}

Theorem Bii is also very close to \eqref{Asymmetric2} since arbitrary small perturbations of row/column square functions $(w > 2)$ dominate in turn row/column maximal functions. The last statement for $x = x_r + x_c$ solves in passing the problem posed in \cite[Section 7.2]{DeJu04} and refines Theorem A. Again, the proof is strongly based on the decomposition of $L_p$ given in Theorem Bi in conjunction with \eqref{AlgDavis} and algebraic atomic descriptions of the involved Hardy spaces. These latter results are apparently new even for classical (commutative) probability spaces. 

\section{\bf Proof of Theorem A}

In this section we prove Theorem A and briefly discuss its optimality. We shall also present its applications in terms of almost uniform convergence. Our first task is to recall the noncommutative Davis decomposition from \cite{JuPe14} and the algebraic atomic description of the involved Hardy spaces.

\subsection{Algebraic Davis decomposition}

Given $p \ge 1$ and a weak-$*$ dense filtration $(\M_n)_{n \ge 1}$ in $(\M,\tau)$, the column martingale Hardy space $\H^c_p(\M)$ is the completion of finite $L_p$-martingales with respect to $$\|x\|_{\H^c_p(\M)} = \Big\| \Big( \sum_{n \ge 1} |d_n(x)|^2 \Big)^{\frac{1}{2}} \Big\|_p \quad \mbox{with} \quad d_1(x) = \E_1(x).$$ The space $h_p^c(\M)$ is also defined in a similar way via the conditioned square function $$\|x\|_{h^c_p(\M)} = \Big\| \Big( \sum_{n \ge 1} \E_{n-1} |d_n(x)|^2 \Big)^{\frac{1}{2}} \Big\|_p \quad \mbox{with} \quad \E_0|d_1(x)|^2 = |\E_1(x)|^2.$$ In what follows, we will say that an operator $x$ affiliated to $\M$ is an \emph{algebraic $h^c_p$-atom} whenever it can be written in the form $x = \sum_{n \ge 1} a_n b_n$, with $a_n$ and $b_n$ satisfying the following  conditions for $1/p = 1/2 + 1/q$:
\begin{itemize}
\item[i)] $\E_n(a_n)=0$, $b_n \in L_q(\M_n)$ for all $n \ge 1$,

\item[ii)] $\displaystyle \summ_n \|a_n\|^2_2 \le 1$ and $\displaystyle \Big\| \Big( \summ_n |b_n|^2 \Big)^{\frac{1}{2}} \Big\|_q \le 1$.
\end{itemize}	
This leads to define the column-atomic Hardy space $h_{p,\mathrm{aa}}^c(\M)$ as the completion in $h_p^c(\M)$ of the space whose unit ball is the absolute convex hull of the family of algebraic $h_p^c$-atoms. This yields $$\|x\|_{h^c_{p,\mathrm{aa}}(\M)} = \inf \Big\{ \sum_{j \ge 1} |\lambda_j| \, : \ x = \sum_{j \ge 1} \lambda_j x_j \mbox{ with } x_j \ h^c_p\mbox{-atoms} \Big\}.$$ The space $h_p^{1_c}(\M)$ was introduced in \cite{JuPe14} to replace the diagonal space $h_p^d(\M)$ in the noncommutative Davis decomposition from \cite{JuMe10,Per0}. The advantage is that we may work with a strictly smaller space. Namely, $h^{1_c}_p(\M)$ is the subspace of all martingale difference sequences in $L_p(\M;\ell^c_1)$. We refer to \cite{JuPe14,Per} for precise definitions ---which we shall not use here--- and focus uniquely in the algebraic atomic description. We call $x$ an \emph{algebraic $h^{1_c}_p$-atom} whenever it can be written as $x = \sum_{n \ge 1} d_n(\alpha_n \beta_n)$, with $\alpha_n$ and $\beta_n$ satisfying:
\begin{itemize}
\item[i)] $\alpha_n \in L_2(\M)$, $\beta_n \in L_q(\M)$ for $n \ge 1$,

\item[ii)] $\displaystyle \summ_n \|\alpha_n\|^2_2 \le 1$ and $\displaystyle \Big\| \Big( \summ_n |\beta_n|^2 \Big)^{\frac{1}{2}} \Big\|_q \le 1$.
\end{itemize}	
As above, we set $$\|x\|_{h^{1_c}_{p,\mathrm{aa}}(\M)} = \inf \Big\{ \sum_{j \ge 1} |\lambda_j| \, : \ x = \sum_{j \ge 1} \lambda_j x_j \mbox{ with } x_j \ h^{1_c}_p\mbox{-atoms} \Big\}$$ and define $h^{1_c}_{p,\mathrm{aa}}(\M)$ accordingly. This is the algebraic Davis decomposition.

\begin{theorem}\label{Davis decomposition}
Let $1\leq p<2$. Then $$\H^c_p(\M) \simeq h^c_p(\M)+h^{1_c}_p(\M) \simeq h^c_{p,\mathrm{aa}}(\M)+h^{1_c}_{p,\mathrm{aa}}(\M).$$ In fact, $h^c_p(\M) \hskip-1pt \simeq \hskip-1pt h^c_{p,\mathrm{aa}}(\M)$ and $h^{1_c}_p(\M) \hskip-1pt \simeq \hskip-1pt h^{1_c}_{p,\mathrm{aa}}(\M)$. \hskip-3pt The same holds for row spaces.
\end{theorem}

\dem 
The argument can be found in \cite[Theorem 5.7]{JuPe14} and \cite[Section 3.6]{Per}. 
\fin

\subsection{Proof of Theorem Ai}

Along this section, we will limit ourselves to prove the column statements since their row analogs are proved similarly. The case $p=2$ follows easily from \cite{Jun02} by interpolation. Indeed, given $0 < \theta < 1$ and according to \cite{JuPa10}, there exists $1 < p_0 < 2 < p_1 < \infty$ and $0 < \eta < 1$ satisfying one of the following isomorphisms $$L_2(\M) = \big[ L_{p_0}(\M), L_{p_1}(\M) \big]_\eta$$ and $$L_2(\M; \ell_\infty^\theta) = \begin{cases} \big[ L_{p_0}(\M;\ell_\infty), L_{p_1}(\M;\ell_\infty^r) \big]_\eta & \mbox{if } \theta < 1/2, \\ \big[ L_{p_0}(\M;\ell_\infty), L_{p_1}(\M;\ell_\infty) \big]_\eta & \mbox{if } \theta = 1/2, \\ \big[ L_{p_0}(\M;\ell_\infty), L_{p_1}(\M;\ell_\infty^c) \big]_\eta & \mbox{if } \theta > 1/2. \end{cases}$$ Thus, Theorem Ai for $p=2$ follows at once from the symmetric and asymmetric Doob inequalities in \cite{Jun02} which we recalled in the Introduction. We may therefore assume in what follows that $1 \le p < 2$ and $1-p/2<\theta<1$ is fixed. According to Theorem \ref{Davis decomposition}, it suffices to prove 
\begin{eqnarray} \label{BasicIneqs1}
\big\| (\E_n(x))_{n \ge 1} \big\|_{L_{p}(\M;\ell^\theta_{\infty})} & \le & c_{p,\theta} \|x\|_{h^c_{p,\mathrm{aa}}(\M)}, \\ \label{BasicIneqs2} \big\| (\E_n(x))_{n \ge 1} \big\|_{L_{p}(\M;\ell^\theta_{\infty})} & \le & c_{p,\theta} \|x\|_{h^{1_c}_{p,\mathrm{aa}}(\M)}.
\end{eqnarray}

\noindent {\bf Proof of (\ref{BasicIneqs1})}. Assume $$x = \sum_{j \ge 1} \lambda_j x_j$$ is an algebraic $h_p^c$-atomic decomposition of $x$ satisfying $\lambda_j \ge 0$ and $\sum_j \lambda_j = 1$ by homogeneity. Recall that $x_j = \sum_{m \ge 1} a_m^j b_m^j$ with $\E_m(a_m^j)=0$ and $b_m^j \in L_q(\M_m)$ where $1/p = 1/2 + 1/q$. Then, $\E_n(x)$ admits the following factorization
\begin{eqnarray*}
\E_n(x) \hskip-5pt & = & \hskip-5pt \sum_{j \ge 1} \lambda_j \sum^n_{k=1} d_k \Big( \sum_{m \ge 1} a^j_m b^j_m \Big) = \sum_{j \ge 1} \sum^n_{k=1} \sum_{m<k} \lambda_j d_k(a^j_m) b^j_m \\ \hskip-5pt & = & \hskip-5pt \Big( \underbrace{\sum_{j \ge 1} \sum_{1 \le m < k \le n} \lambda_j^{\frac{1}{2}}d_k(a^j_m) \otimes e_{1(j,m)}}_{A_n} \Big) \Big( \underbrace{\sum_{j \ge 1} \sum_{m \ge 1} \lambda_j^{\frac{1}{2}} b^j_m \otimes e_{(j,m)1}}_{B} \Big).
\end{eqnarray*}
Using again $\E_m(a^j_m)=0$ we observe that $$A_n = \sum^n_{k=1} \Big( \sum_{j \ge 1} \sum_{m=1}^{k-1} \lambda_j^{\frac{1}{2}} d_k(a^j_m) \otimes e_{1(j,m)} \Big) = \E_n\otimes id \Big( \sum_{j \ge 1} \sum_{m \ge 1} \lambda_j^{\frac{1}{2}} a^j_m \otimes e_{1(j,m)} \Big) = \widehat{\E}_n(A)$$ where $\widehat{\E}_n = \E_n \otimes id_{\mathcal{B}(\ell_2)}$. By the definition of algebraic atoms $$\|A\|_2 = \Big( \sum_{j,m \ge 1} \big\| \lambda_j^{\frac{1}{2}}a^j_m \big\|^2_2 \Big)^{\frac12} = \Big( \sum_{j \ge 1} \lambda_j \sum_{m \ge 1} \|a^j_m\|^2_2 \Big)^{\frac12} \le 1.$$ Letting $\frac{1}{s} = \frac12 - \frac{1-\theta}{p} > 0$, polar decomposition yields $A = A_{\mathrm{row}} A_{\mathrm{mat}}$ where 
\begin{itemize}
\item[a)] $A_{\mathrm{row}} \in L_{\frac{p}{1-\theta}}(\M \bar\otimes \mathcal{B}(\ell_2))$ is a row matrix with $\|A_{\mathrm{row}}\|_{\frac{p}{1-\theta}} \le 1$,

\item[b)] $A_{\mathrm{mat}} \in L_s(\M \bar\otimes \mathcal{B}(\ell_2))$ is a full matrix satisfying $\|A_{\mathrm{mat}}\|_s \le 1$.
\end{itemize}
According to \cite[Proposition 2.8]{Jun02}, for each $n \ge 1$ there is an isometric right $\M_n\bar{\otimes}\mathcal{B}(\ell_2)$-module map $\widehat{u}_{n}: \M \bar\otimes \mathcal{B}(\ell_2) \to C(\M_n\bar{\otimes}\mathcal{B}(\ell_2))$ whose image is the space of columns with entries in $\M_n\bar{\otimes}\mathcal{B}(\ell_2)$ and such that $$\widehat{\E}_n(A) = \widehat{u}_{n}(A_{\mathrm{row}}^*)^* \hskip1pt \widehat{u}_{n} ( A_{\mathrm{mat}}^{}).$$ On the other hand, by the symmetric Doob maximal inequality \cite{Jun02} in the amplified space $L_{\frac{p}{2(1-\theta)}}(\M\bar{\otimes}\mathcal{B}(\ell_2))$, we find $\alpha \in L_{\frac p{1-\theta}}(\M \bar\otimes \mathcal{B}(\ell_2))$ and $\rho_n \in \M \bar\otimes \mathcal{B}(\ell_2)$ which satisfy the following relations for $n \ge 1$
$$\widehat{\E}_n(A_{\mathrm{row}}^{} A_{\mathrm{row}}^*) = \widehat{u}_n(A_{\mathrm{row}}^*)^* \widehat{u}_n(A_{\mathrm{row}}^*) = \alpha^* \rho_n^* \rho_n \alpha,$$ $$\| \alpha\|_{\frac p{1-\theta}} \Big( \sup_{n \ge 1} \|\rho_n\|_\infty \Big) \le C_{\frac p{2(1-\theta)}} \big\| A_{\mathrm{row}}^{} A_{\mathrm{row}}^* \big\|_{\frac p{2(1-\theta)}}^{\frac12} = C_{\frac p{2(1-\theta)}} \|A_{\mathrm{row}}\|_{\frac{p}{1-\theta}} \le c_{p,\theta}^1.$$ Similarly, we may find $\beta \in L_s(\M \bar\otimes \mathcal{B}(\ell_2))$ and $\gamma_n \in \M \bar\otimes \mathcal{B}(\ell_2)$ for any $n \ge 1$ with $$\widehat{u}_n(A_{\mathrm{mat}})^* \widehat{u}_n(A_{\mathrm{mat}}^{}) = \beta^* \gamma_n^* \gamma_n \beta \quad \mbox{and} \quad \|\beta\|_s \Big( \sup_{n \ge 1} \|\gamma_n\|_\infty \Big) \le c_{p,\theta}^2.$$ According to polar decomposition and the factorizations found so far, it is not difficult to construct contractions $\xi_n, \psi_n \in \M \bar\otimes \mathcal{B}(\ell_2)$ so that $\E_n(x)$ may be rewritten as follows 
\begin{eqnarray*}
\E_n(x) & = & \widehat{\E}_n(A) B \\ [3pt] & = & \widehat{u}_{n}(A_{\mathrm{row}}^*)^* \hskip1pt \widehat{u}_{n} ( A_{\mathrm{mat}}^{}) B \ = \ \big( \alpha^* \rho_n^* \xi_n^*\big) \big( \psi_n \gamma_n \beta \big) B \\ & = & \underbrace{(\alpha^* \alpha)^{\frac12}}_{a} \Big( \underbrace{(\alpha^* \alpha)^{-\frac12} \alpha^* \big( \rho_n^* \xi_n^* \psi_n \gamma_n \big) \beta B (B^* \beta^* \beta B)^{-\frac12}}_{w_n} \Big) \underbrace{(B^* \beta^* \beta B)^{\frac12}}_{b}.
\end{eqnarray*}
We claim that $(a, w_n, b) \in L_{\frac{p}{(1-\theta)}}(\M) \times L_\infty(M) \times L_{\frac{p}{\theta}}(\M)$ and $$\|a\|_{\frac{p}{1-\theta}} \Big( \sup_{n \ge 1} \|w_n\|_\infty \Big) \|b\|_{\frac{p}{\theta}} \le c_{p,\theta}^1 c_{p,\theta}^2 = c_{p,\theta}.$$ This implies \eqref{BasicIneqs1}. The fact that $a,b,w_n$ are affiliated with $\M$ and not with the amplified algebra $\M \bar\otimes \mathcal{B}(\ell_2)$ boils down to the observation that $\alpha^*$ is a row matrix and $B$ a column matrix. Note that $\alpha^*$ is a row matrix because the same holds for $\widehat{u}_n(A_{\mathrm{row}}^*)^*$ since its adjoint is a column of columns. On the  other hand, since $\xi_n^*$ and $\psi_n$ as well as $(\alpha^* \alpha)^{-1/2} \alpha^*$ and $\beta B (B^* \beta^* \beta B)^{-1/2}$ are contractions, we conclude that $$\|a\|_{\frac{p}{1-\theta}} \Big( \sup_{n \ge 1} \|w_n\|_\infty \Big) \|b\|_{\frac{p}{\theta}} \le \|\alpha\|_{\frac{p}{1-\theta}} \Big( \sup_{n \ge 1} \|\rho_n\|_\infty \|\gamma_n\|_\infty \Big) \|\beta\|_s \|B\|_q \le c_{p,\theta}^1 c_{p,\theta}^2.$$ Indeed, $\frac{\theta}{p} = \frac1s + \frac1q$ and $\displaystyle \|B\|_q = \Big\| \sum_{j,m \ge 1} \lambda_j |b^j_m|^2 \Big\|^{\frac12}_{\frac{q}{2}} \le \Big( \sum_{j \ge 1} \lambda_j \Big\| \sum_{m \ge 1} |b^j_m|^2 \Big\|_{\frac{q}{2}} \Big)^{\frac12} \le 1$.

\noindent {\bf Proof of (\ref{BasicIneqs2})}. Assume $$x = \sum_{j \ge 1} \lambda_j x_j$$ is an algebraic $h_p^{1_c}$-atomic decomposition of $x$ satisfying $\lambda_j \ge 0$ and $\sum_j \lambda_j = 1$ by homogeneity. Recall $x_j = \sum_{m \ge 1} d_m(\alpha_m^j \beta_m^j)$ with $\alpha_m^j \in L_2(\M)$ and $\beta_m^j \in L_q(\M)$ where $1/p = 1/2 + 1/q$. Then, $\E_n(x)$ may be written as
\begin{eqnarray*}
\E_n(x) & = & \sum^n_{k=1}d_k(x)=\sum_{j\ge1} \lambda_j\sum^n_{k=1}d_k(\alpha^j_k \beta^j_k) \\
& = & \sum_{j \ge 1} \lambda_j \sum^n_{k=1} \E_k(\alpha^j_k \beta^j_k) - \sum_{j\ge1} \lambda_j \sum^n_{k=1} \E_{k-1}(\alpha^j_k \beta^j_k) \ = \ X_n - Y_n.
\end{eqnarray*}
By the (quasi)-triangle inequality in $L_p(\M;\ell_\infty^\theta)$, it suffices to estimate the norms of $(X_n)_{n \ge 1}$ and $(Y_n)_{n \ge 1}$ separately. Since both are similar, we shall only justify the one for $X_n$'s. To that end we emulate the argument for $h_p^c$-atoms, so that we aim to express $X_n$ in the form $\mathbb{E}_n(A)B$ for some operators $A,B$ affiliated to $\M \bar\otimes \mathcal{B}(\ell_2)$ and some conditional expectations $\mathbb{E}_n$. The cancelation of $h_p^c$-atoms allowed us to take $\mathbb{E}_n = \E_n \otimes id_{\mathcal{B}(\ell_2)}$ above. Our choice this time will be different. Before that we apply \cite[Proposition 2.8]{Jun02} to factorize $$\E_{k}(\alpha^j_k \beta^j_k) = u_{k}(\alpha^{j*}_k)^* u_{k}(\beta^j_k)$$ where $u_{k}: \M \to C(\M_{k})$ is an isometric right $\M_{k}$-module map. Hence $$X_n = \Big( \underbrace{\sum_{j \ge 1} \sum_{k \le n} \lambda_j^{\frac12} (u_{k}(\alpha^{j*}_k))^*\otimes e_{1j}\otimes e_{1k}}_{A_n} \Big) \Big( \underbrace{\sum_{j \ge 1} \sum_{k \ge 1} \lambda_j^{\frac12} u_{k}(\beta^j_k) \otimes e_{j1} \otimes e_{k1}}_B \Big).$$ Let us take $\mathbb{E}_n = id_\M \otimes id_{\mathcal{B}(\ell_2)} \otimes id_{\mathcal{B}(\ell_2)} \otimes \mathsf{E}_n$ where $$\mathsf{E}_n \big( (m_{jk})_{j,k \ge 1} \big) = (m_{jk})_{1 \le j,k \le n} \oplus (m_{kk})_{k>n}$$ is a unital conditional expectation in $\mathcal{B}(\ell_2)$. Of course, this gives $A_n = \mathbb{E}_n(A)$ as desired. Once this is clarified, the estimate for the $L_p(\M; \ell_\infty^\theta)$-norm of $(X_n)_{n \ge 1}$ can be deduced following the same argument we used for $h_p^c$-atoms above as long as we can prove that $\|A\|_2$ and $\|B\|_q$ are finite. We have 
\begin{eqnarray*}
\|A\|_2 \hskip-5pt & = & \hskip-5pt \Big( \sum_{j \ge 1} \lambda_j \sum_{k \ge 1} \tau \big( \E_{k}(|\alpha^j_k|^2) \big) \Big)^{\frac12} = \hskip1pt \Big( \sum_{j \ge 1} \lambda_j \sum_{k \ge 1} \|\alpha^j_k\|_2^2 \Big)^{\frac12} \le 1, \\ \|B\|_q \hskip-5pt & = & \hskip-5pt \Big\| \sum_{j \ge 1} \lambda_j \sum_{k\geq1} \E_{k}(|\beta^j_k|^2) \Big\|_{\frac{q}{2}}^{\frac12} \le \Big( C_{\frac{q}{2}} \sum_{j \ge 1} \lambda_j \Big\| \sum_{k \ge 1} |\beta^j_k|^2 \Big\|_{\frac{q}{2}} \Big)^{\frac12} \le c_q.
\end{eqnarray*}
The bound of $B$ follows from the dual Doob inequalities \cite{Jun02} since $1 \le q/2 < \infty$. \fin

\subsection{Proof of Theorem Aii}

As above, it suffices to consider the column spaces and we begin with the case $p=2$. By the definition of ${\Lambda_{2,\infty}(\M;\ell^c_{\infty})}$, it can be easily checked that
\begin{eqnarray*}
\big\| (\E_n(x))_{n \ge 1} \big\|_{\Lambda_{2,\infty}(\M;\ell^c_{\infty})} \hskip-6pt & = & \hskip-6pt \big\| (|\E_n(x)|^2)_{n \ge 1} \big\|^{\frac12}_{\Lambda_{1,\infty}(\M;\ell_{\infty})} \\ \hskip-6pt & \le & \hskip-6pt \big\| (\E_n(|x|^2))_{n \ge 1} \big\|^{\frac12}_{\Lambda_{1,\infty}(\M;\ell_{\infty})} \le \| \hskip1pt |x|^2 \|_{L_1(\M)}^{\frac12} = \|x\|_{L_2(\M)}. 
\end{eqnarray*}
Here we have used Kadison-Schwarz inequality and Cuculescu weak type estimate \cite{Cu}, which holds with constant $1$. This proves the result for $p=2$ since we have $\H_2^c(\M) = L_2(\M)$. Let us now assume that $1 \le p < 2$. By the algebraic Davis decomposition in Theorem \ref{Davis decomposition}, it suffices to prove
\begin{eqnarray}
\label{BasicIneqs3}
\big\| (\E_n(x))_{n \ge 1} \big\|_{\Lambda_{p,\infty}(\M;\ell^c_{\infty})} & \le & c_p \|x\|_{h^c_{p,\mathrm{aa}}(\M)} \\ 
\label{BasicIneqs4}
\big\| (\E_n(x))_{n \ge 1} \big\|_{\Lambda_{p,\infty}(\M; \ell^c_{\infty})} & \le & c_p \|x\|_{h^{1_c}_{p,\mathrm{aa}}(\M)}.
\end{eqnarray}

\noindent {\bf Proof of (\ref{BasicIneqs3})}. Assume by homogeneity that $$\|x\|_{h_{p,\mathrm{aa}}^c(\M)} < 1$$ and follow the proof of Theorem Ai to factorize $\E_n(x) = \widehat{\E}_n(A) B$ with $$\max \Big\{ \|A\|_{L_2(\M \bar\otimes \mathcal{B}(\ell_2))}, \|B\|_{L_q(\M \bar\otimes \mathcal{B}(\ell_2))} \Big\} < 1.$$ According to Theorem Aii for $p=2$ (already justified with $c_2=1$) we obtain $$\big\| ( \widehat{\E}_n(A) )_{n \ge 1} \big\|_{\Lambda_{2,\infty}(\M \bar\otimes \mathcal{B}(\ell_2))} \le \|A\|_{L_2(\M \bar\otimes \mathcal{B}(\ell_2))} < 1.$$ We are now ready to justify \eqref{BasicIneqs3}. Indeed, given $\lambda>0$ set in what follows $\lambda_1 = \lambda^{p/2}$ and $\lambda_2 = \lambda^{p/q}$. According to the definition of the weak space $\Lambda_{2,\infty}(\M \bar\otimes \mathcal{B}(\ell_2); \ell_\infty^c)$ there must exist a projection $e_{\lambda_1} \in \M \bar\otimes \mathcal{B}(\ell_2)$ satisfying $$\big| \widehat{\E}_n(A) e_{\lambda_1} \big| \le \lambda_1 \hskip5pt \mbox{and} \hskip5pt \lambda_1 \big( \widehat{\tau}(\1 - e_{\lambda_1}) \big)^{\frac{1}{2}} \le (1+\delta) \big\| (\widehat{\E}_n(A))_{n \ge 1} \big\|_{\Lambda_{2,\infty}(\M \bar\otimes \mathcal{B}(\ell_2); \ell^c_{\infty})} < 1,$$ where $\widehat{\tau} = \tau \otimes \mathrm{tr}$. In addition, $B$ is a column so that $|B| \in L_q(\M)$. This means that the spectral projection $f_{\lambda_2} = \chi_{[0,\lambda_2]}(|B|)$ belongs to $\M$. Moreover, by Chebyshev inequality we also find that the following inequalities hold $$|B| f_{\lambda_2} \le \lambda_2 \quad \mbox{and} \quad \lambda_2 ( \tau(\1 - f_{\lambda_2}))^{\frac1q} \le \|B\|_q < 1.$$ Then we construct the following projection in $\M$ $$\Pi_{\lambda} = \Big( \1 - \mathrm{supp} \big| (\1-e_{\lambda_1})B \big| \Big) \wedge f_{\lambda_2}.$$ Observe that $(\1 - e_{\lambda_1}) B \Pi_{\lambda}=0$, which yields in turn $$\widehat{\E}_n(A) B \Pi_{\lambda} = \widehat{\E}_n(A) e_{\lambda_1} B \Pi_{\lambda} = \widehat{\E}_n(A) e_{\lambda_1} B f_{\lambda_2} \Pi_{\lambda} \Rightarrow \big\| \widehat{\E}_n(A) B \Pi_{\lambda} \big\|_{\infty} \le \lambda_1\lambda_2 = \lambda.$$ Therefore, by the definition of $\Lambda_{p,\infty}(\M; \ell_\infty^c)$ it suffices to estimate $\lambda (\tau(\1-\Pi_{\lambda}))^{\frac{1}{p}}$
\begin{eqnarray*}
\lambda \big( \tau(\1 - \Pi_\lambda) \big)^{\frac{1}{p}} \hskip-3pt & \le & \hskip-3pt \lambda \Big( \tau \big( \mathrm{supp}|(\1-e_{\lambda_1})B| \big) + \tau \big( \1 -f_{\lambda_2} \big) \Big)^{\frac{1}{p}} \\ \hskip-3pt & = & \hskip-3pt \lambda \Big( \widehat{\tau} \big( \mathrm{supp}|B^*(\1-e_{\lambda_1})| \big) + \tau \big( \1 - f_{\lambda_2} \big) \Big)^{\frac{1}{p}} \\ \hskip-3pt & \le & \hskip-3pt \lambda \Big( \widehat{\tau}(\1 - e_{\lambda_1}) + \tau (\1-f_{\lambda_2}) \Big)^{\frac{1}{p}} \ < \ \lambda \Big( \lambda_1^{-2} + \lambda_2^{-q} \Big)^{\frac{1}{p}} \ = \ 2^{\frac1p}.
\end{eqnarray*}

\noindent {\bf Proof of (\ref{BasicIneqs4})}. Assume by homogeneity that $$\|x\|_{h_{p,\mathrm{aa}}^{1_c}(\M)} < 1$$ and follow the proof of Theorem Ai to write $\E_n(x) = X_n - Y_n$, where both $X_n$ and $Y_n$ are of the form $\mathbb{E}_n(A)B$ for certain rows $A$ and columns $B$ satisfying the same estimates above $$\max \Big\{ \|A\|_{L_2(\M \bar\otimes \mathcal{B}(\ell_2))}, \|B\|_{L_q(\M \bar\otimes \mathcal{B}(\ell_2))} \Big\} < 1.$$ According to Theorem Aii for $p=2$ we obtain $$\big\| ( \mathbb{E}_n(A) )_{n \ge 1} \big\|_{\Lambda_{2,\infty}(\M \bar\otimes \mathcal{B}(\ell_2))} \le \|A\|_{L_2(\M \bar\otimes \mathcal{B}(\ell_2))} < 1.$$ Then, our argument above for \eqref{BasicIneqs3} applies and yields the following inequalities $$\max \Big\{ \big\| (X_n)_{n \ge 1} \big\|_{\Lambda_{p,\infty}(\M;\ell_\infty^c)}, \big\| (Y_n)_{n \ge 1} \big\|_{\Lambda_{p,\infty}(\M;\ell_\infty^c)} \Big\} \le c_p.$$ The desired result follows from the quasi-triangle inequality in $\Lambda_{p,\infty}(\M;\ell_\infty^c)$. \fin

\begin{remark} 
\emph{The idea behind the proof of \eqref{BasicIneqs3} is a H\"older type inequality for the Cuculescu spaces $\Lambda_{p,\infty}(\M;\ell_\infty)$. More precisely, given $0 < p,q\le \infty$ such that $1/r=1/p+1/q$, we find that $$\Lambda_{p,\infty}(\M;\ell^c_{\infty})L_q(\M)\subset\Lambda_{r,\infty}(\M;\ell^c_{\infty}).$$
In other words, the following inequality holds $$\big\| (x_nb)_{n \ge 1} \big\|_{\Lambda_{r,\infty}(\M;\ell^c_{\infty})} \le 2^{\frac1r} \big\| (x_n)_{n \ge 1} \big\|_{\Lambda_{p,\infty}(\M;\ell^c_{\infty})} \|b\|_{L_q(\M)}.$$ The proof can be reconstructed from our proof of \eqref{BasicIneqs3}, but the argument there is a bit more involved since our operators $x_n$ and $b$ live in the matrix amplified algebra $\M \bar\otimes \mathcal{B}(\ell_2)$ although their product does not. This forces us to be a bit more careful.}
\end{remark}

\subsection{Conclusions}

We conclude this section with a little discussion on $\tau$-almost uniform convergence and the optimality of Theorem A. Let us precise our definition of almost uniform convergence given in the Introduction. A sequence $(x_n)_{n \ge 1}$ of $\tau$-measurable operators converges to $0$ \emph{$\tau$-almost uniformly from the right} when there is a sequence of projections $(p_k)_{k \ge 1}$ in $\M$ satisfying $\lim_k \tau( \1 - p_k) = 0$ and $\lim_n \|x_n p_k\|_\infty = 0$ for all $k \ge 1$. Similarly, $(x_n)_{n \ge 1}$ converges to $0$ \emph{$\tau$-almost uniformly from the left} when $\lim_n \|p_k x_n\|_\infty = 0$ instead.

\begin{corollary}
Given $1 \le p \le 2$ and $x \in \H_p^c(\M)$, the sequence $\E_n(x)$ converges $\tau$-almost uniformly from the right to $x$. Similarly, when $x \in \H_p^r(\M)$ the $\tau$-a.u. convergence holds from the left.
\end{corollary}

\dem Recall from Theorem Aii that 
\begin{eqnarray*}
\lefteqn{\big\| (\E_n(x) - \E_m(x))_{n \ge m} \big\|_{\Lambda_{p,\infty}(\M,\ell_\infty^c)}} \\ & = & \big\| (\E_n(x - \E_m(x)))_{n \ge m} \big\|_{\Lambda_{p,\infty}(\M,\ell_\infty^c)} \ \le \ \big\| x - \E_m(x) \big\|_{\H_p^c(\M)} \ \to \ 0
\end{eqnarray*}
as $m \to \infty$. Combining this with the proof of \cite[Proposition 5.1]{DeJu04} we obtain the desired result. The row case is justified similarly. This completes the proof. \fin

\begin{remark} \label{Optimalidad}
\emph{According to \cite{JX0}, the symmetric estimate $\H_1(\M) \to L_1(\M; \ell_\infty)$ fails and our restrictions $\theta < p/2$ in the row case and $\theta > 1 - p/2$ in the column case become necessary for $p=1$. In addition, since $L_2(\M) = \H_2^r(\M) = \H_2^c(\M)$ the negative results in \cite{DeJu04} for $p < 2$ indicate that we may not expect a better result for $p=2$. These considerations lead us to conjecture that Theorem A is best possible in our restrictions for the parameter $0 \le \theta \le1$.}
\end{remark}

\section{\bf Proof of Theorem B}

In this section we prove Theorem B. This requires to introduce a family of Hardy spaces, apparently new even in classical/commutative martingale $L_p$ theory. As a crucial point in our approach, we shall also investigate their dual spaces.

\subsection{New Hardy spaces}

The family of Hardy spaces to be introduced below is motivated by an elementary observation. Namely, that the norms in $h_{p,\mathrm{aa}}^c(\M)$ and $h_{p,\mathrm{aa}}^{1_c}(\M)$ can be simplified as follows.

\begin{lemma} \label{Simplification}
Given $1 \le p \le 2$ with $1/p = 1/2 + 1/q$, we find
\begin{eqnarray*}
\|x\|_{h^c_{p,\mathrm{aa}}(\M)} & = & \inf_{\begin{subarray}{c} x=\sum_n a_nb_n \\ \E_n(a_n)=0, b_n \in L_q(\M_n)\end{subarray}} \hskip1pt \Big( \sum_{n \ge 1} \|a_n\|^2_2 \hskip1pt \Big)^{\frac12} \hskip2pt \Big\| \Big( \sum_{n \ge 1} \hskip1pt |b_n|^2 \hskip1pt \Big)^{\frac12} \Big\|_q, \\ \|x\|_{h^{1_c}_{p,\mathrm{aa}}(\M)} & = & \inf_{\begin{subarray}{c} x=\sum_n d_n(\alpha_n \beta_n) \\ \alpha_n \in L_2(\M), \beta_n \in L_q(\M) \end{subarray}} \Big( \sum_{n \ge 1} \|\alpha_n\|^2_2 \Big)^{\frac12} \hskip1pt \Big\| \Big( \sum_{n \ge 1} |\beta_n|^2 \Big)^{\frac12} \Big\|_q.
\end{eqnarray*}
In other words, any sum of atoms $\sum_j \lambda_j x_j$ may be rewritten as a single algebraic $h_p^c$-atom or $h_p^{1_c}$-atom accordingly. Similar simplifications apply for the  row spaces. 
\end{lemma}

\dem The proof is very similar in all cases, let us justify it rigorously for the space $h_{p,\mathrm{aa}}^c(\M)$. The quantity on the right hand side is clearly not smaller than the one on the left hand side. Conversely, given $\delta > 0$ consider a decomposition $x = \sum_j \lambda_j x_j$ into algebraic $h_p^c$-atoms satisfying $\lambda_j > 0$ and $\sum_j\lambda_j \le (1+\delta)\|x\|_{h_{p,\mathrm{aa}}^c(\M)}$. We know $x_j = \sum_m a_m^j b_m^j$ with $\E_m(a_m^j)=0$ and $b_m^j \in L_q(\M_m)$. Set $$a_m = \sum_{j \ge 1} \lambda_j a^j_m b^j_m b_m^{-1} \quad \mbox{where} \quad b_m = \Big( \sum_{j \ge 1} \lambda_j |b^j_m|^2 \Big)^{\frac12}.$$
Recall that $\E_n(a_n)=0$ and $b_n \in L_q(\M_n)$. Thus, it suffices to prove
\begin{itemize}
\item[i)] $x = \displaystyle \summ_n a_n b_n$, 

\item[ii)] $\displaystyle \Big( \summ_n \|a_n\|_2^2 \Big)^{\frac12} \Big\| \Big( \summ_n |b_n|^2 \Big)^\frac12 \Big\|_q \le \summ_j \lambda_j$,
\end{itemize}
and conclude by letting $\delta \to 0^+$. The first identity requires to justify Fubini in  
\begin{eqnarray*}
x = \sum_{j \ge 1} \lambda_j x_j \hskip-5pt & = & \hskip-5pt h_{p,\mathrm{aa}}^c-\lim_{J \to \infty} \sum_{j < J} \sum_{m \ge 1} \lambda_j a_m^j b_m^j \\ \hskip-5pt & = & \hskip-5pt  h_{p,\mathrm{aa}}^c-\lim_{J \to \infty} \sum_{m \ge 1} \sum_{j < J} \lambda_j a_m^j b_m^j = \sum_{m \ge 1} a_m b_m.
\end{eqnarray*}
The first limit holds since the partial sums $\sum_{j < J} \lambda_j \sum_m a_m^j b_m^j$ are clearly Cauchy in $h_{p,\mathrm{aa}}^c(\M)$. The second limit requires to show that the $h_{p,\mathrm{aa}}^c$-norms of the following sums converge to $0$ as $J \to \infty$ $$\sum_{m \ge 1} \sum_{j \ge J} \lambda_j a_m^j b_m^j = \sum_{m \ge 1} a_m(J) b_m(J)$$ where $a_m(J) = \sum_{j \ge J} \lambda_j a_m^j b_m^j b_m(J)^{-1}$ and $b_m(J)^2 = \sum_{j \ge J} \lambda_j |b_m^j|^2$. This is an algebraic $h_p^c$-atom since $\E_m(a_m(J))=0$ and $b_m(J) \in L_q(\M_m)$. In particular, we immediately deduce the following estimate $$\Big\| \sum_{m \ge 1} \sum_{j \ge J} \lambda_j a_m^j b_m^j  \Big\|_{h_{p,\mathrm{aa}}^c(\M)} \le \Big( \sum_{m \ge 1} \|a_m(J)\|_2^2 \Big)^\frac12 \hskip2pt \Big\| \Big( \sum_{m \ge 1} |b_m(J)|^2 \Big)^\frac12 \Big\|_q.$$ We may estimate the first term on the right hand side as follows
\begin{eqnarray*}
\sum_{m \ge 1} \|a_m(J)\|_2^2 \hskip-5pt & \le & \hskip-5pt \sum_{m \ge 1} \Big\| \Big( \sum_{j \ge J} \lambda_j |a_m^{j*}|^2 \Big)^\frac12 \Big\|_2^2 \Big\| \Big( \sum_{j \ge J} \lambda_j |b_m^j b_m(J)^{-1}|^2 \Big)^\frac12 \Big\|_\infty^2 \\ \hskip-5pt & \le & \hskip-5pt \sum_{m \ge 1} \Big\| \Big( \sum_{j \ge J} \lambda_j |a_m^{j*}|^2 \Big)^\frac12 \Big\|_2^2 = \sum_{m \ge 1} \sum_{j \ge J} \lambda_j \|a_m^j\|_2^2 \le \sum_{j \ge J} \lambda_j.
\end{eqnarray*}
The second term uses the triangle inequality in $L_{q/2}(\M)$ since $q \ge 2$
\begin{eqnarray*}
\Big\| \Big( \sum_{m \ge 1} |b_m(J)|^2 \Big)^\frac12 \Big\|_q^2 \hskip-5pt & = & \hskip-5pt \Big\| \sum_{m \ge 1} \sum_{j \ge J} \lambda_j |b_m^j|^2 \Big\|_{\frac{q}{2}} \le \sum_{j \ge J} \lambda_j \Big\| \Big( \sum_{m \ge 1} |b_m^j|^2 \Big)^\frac12 \Big\|_q^2 \le \sum_{j \ge J} \lambda_j.
\end{eqnarray*}
Altogether we obtain $$\lim_{J \to \infty} \Big\| \sum_{m \ge 1} \sum_{j \ge J} \lambda_j a_m^j b_m^j  \Big\|_{h_{p,\mathrm{aa}}(\M)} \le \lim_{J \to \infty} \sum_{j \ge J} \lambda_j = 0,$$ which completes the proof of claim i). Claim ii) follows from above for $J=1$. The assertion for $h_{p,\mathrm{aa}}^{1_c}(\M)$ is very similar. Indeed, given $x = \sum_j \lambda_j \sum_m d_m (\alpha_m^j \beta_m^j)$ pick $$\alpha_m = \summ_j \lambda_j \alpha_m^j \beta_m^j \beta_m^{-1} \quad \mbox{with} \quad \beta_m = \Big( \summ_j \lambda_j |\beta_m^j|^2 \Big)^\frac12.$$ The exact same argument yields a suitable decomposition $x = \sum_m d_m(\alpha_m \beta_m)$. \fin

We are now ready to generalize the family of algebraic atomic Hardy spaces. Let $L_0(\M,\tau)$ stand for the space of $\tau$-measurable operators. Let $1 \le p \le 2$ and $s \ge 2$ so that $1/p=1/w+1/s$. Then we define 
\begin{eqnarray*}
h^c_{pw}(\M) & = & \Big\{ x \in L_0(\M,\tau) \, : \ \|x\|_{h^c_{pw}(\M)} < \infty \Big\}, \\ h^{1_c}_{pw}(\M) & = & \Big\{ x \in L_0(\M,\tau) \, : \ \|x\|_{h^{1_c}_{pw}(\M)} <\infty \Big\}, 
\end{eqnarray*}
where
\begin{eqnarray*}
\|x\|_{h^c_{pw}(\M)} & = & \inf_{\begin{subarray}{c} x = \sum_n a_n b_n \\ \E_n(a_n)=0, \ b_n \in L_s(\M_n) \end{subarray}} \Big\| \sum_{n \ge 1} \hskip1pt a_n \hskip1pt \otimes e_{1n} \Big\|_w \Big\| \sum_{n \ge 1} \hskip1pt b_n \otimes e_{n1} \Big\|_s, \\ \|x\|_{h^{1_c}_{pw}(\M)} & = & \inf_{\begin{subarray}{c} x = \sum_n d_n(\alpha_n \beta_n) \\ \alpha_n \in L_w(\M), \beta_n \in L_s(\M) \end{subarray}} \Big\| \sum_{n \ge 1} \alpha_n \otimes e_{1n} \Big\|_w \Big\| \sum_{n \ge 1} \beta_n \otimes e_{n1} \Big\|_s.
\end{eqnarray*}
The analog families of row Hardy spaces are defined by taking adjoints as usual. 

\begin{remark} \label{w=2}
\emph{According to Theorem \ref{Davis decomposition} and Lemma \ref{Simplification}, we get the isomorphisms} 
\begin{eqnarray*}
h_{p2}^c(\M) \hskip-3pt & = & \hskip-3pt h_{p,\mathrm{aa}}^c(\M) \ \simeq \ h_p^c(\M), \\ h_{p2}^{1_c}(\M) \hskip-3pt & = & \hskip-3pt h_{p,\mathrm{aa}}^{1_c}(\M) \ \simeq \ h_p^{1_c}(\M).
\end{eqnarray*}
\emph{Of course we could have allowed $s \ge p$ by imposing $w \ge 2$. The most interesting spaces for this paper will be those satisfying $w,s \ge 2$, although those with $w < 2$ will also be instrumental for our purposes.}  
\end{remark}
\begin{lemma}\label{lem:norm}
The following holds:
\begin{itemize}
\item[i)] If $w \ge 2$, $\|\cdot\|_{h^c_{pw}(\M)}$ is a norm.

\item[ii)] If $w < 2$, $\|\cdot\|_{h^c_{pw}(\M)}$ is a $\frac{w}{2}$-norm. 
\end{itemize}
The same is true for $\|\cdot\|_{h^{1_c}_{pw}(\M)}$ and the corresponding row analogues.
\end{lemma}

\dem Homogeneity is straightforward, while the positive definiteness follows from the simple fact that $h_{pw}^c(\M)$ embeds in $L_p(\M)$. The same embedding holds for $h_{pw}^{1_c}(\M)$ ---and so positive definiteness--- although this will require a more involved argument in Lemma \ref{lem:hpw to lp} below. It remains to justify the triangle inequality for $w \ge 2$ and its $\frac{w}{2}$-analogue for $w < 2$. Since the argument is similar for $h_{pw}^c(\M)$ and $h_{pw}^{1_c}(\M)$ ---see the end of the proof of Lemma \ref{Simplification} for a similar arguing--- we shall only consider the space $h_{pw}^c(\M)$. Given $x_1, x_2 \in h_{pw}^c(\M)$ and $\delta > 0$, we may write $x_j = \sum_m a_m^j b_m^j$ with $\E_m(a^j_m) = 0$ and $b^j_m \in L_s(\M_n)$ such that $$\Big\| \sum_{m \ge 1} a^j_m \otimes e_{1m} \Big\|_w \Big\| \sum_{m \ge 1} b^j_m \otimes e_{m1} \Big\|_s \le (1+\delta) \|x_j\|_{h^c_{pw}(\M)}.$$
Moreover, by renormalization we may assume $$\Big\| \sum_{m \ge 1} a_m^j \otimes e_{1m} \Big\|_w = \Big\| \sum_{m \ge 1} b^j_m \otimes e_{m1} \Big\|_s.$$ Set $a_m = \sum^2_{j=1} a^j_m b^j_m b_m^{-1}$ and $b_m=(|b^1_m|^2+|b^2_m|^2)^{\frac12}$. This allows us to write $x_1+x_2 = \sum_{m \ge 1} a_m b_m$ with $\E_m(a_m)=0$ and $b_m \in L_s(\M_m)$. Assume now that $w < 2$, our considerations so far yield $$\big\| x_1 + x_2 \big\|_{h_{pw}^c(\M)} \le \Big\| \sum_{m \ge 1} a_m \otimes e_{1m} \Big\|_w \Big\| \sum_{m \ge 1} b_m \otimes e_{m1} \Big\|_s.$$ Therefore, assertion ii) will follow by letting $\delta \to 0^+$ if we can prove $$\max \Big\{ \Big\| \sum_{m \ge 1} a_m \otimes e_{1m} \Big\|_w^w, \Big\| \sum_{m \ge 1} b_m \otimes e_{m1} \Big\|_s^w  \Big\} \le (1+\delta)^{\frac{w}{2}} \sum_{j=1}^2 \|x_j\|_{h^c_{pw}(\M)}^{\frac{w}{2}}$$ or equivalently 
\begin{eqnarray*}
A \ = \ \Big\| \sum_{m \ge 1} a_m \otimes e_{1m} \Big\|_w & \le & \Big( \sum_{j=1}^2 \Big\| \sum_{m \ge 1} a_m^j \otimes e_{1m} \Big\|_w^w \Big)^{\frac{1}{w}}, \\ B \ = \ \Big\| \sum_{m \ge 1} \hskip1pt b_m \otimes e_{m1} \Big\|_s \hskip2pt & \le & \Big( \sum_{j=1}^2 \Big\| \sum_{m \ge 1} \hskip1pt b_m^j \otimes e_{m1} \Big\|_s^w \Big)^{\frac{1}{w}}.
\end{eqnarray*} 
To prove the first estimate we note that 
\begin{eqnarray*}
A & = & \Big\| \sum_{j=1}^2 \Big( \sum_{m \ge 1} a_m^j \otimes e_{1m} \Big) \Big( \sum_{m \ge 1} b_m^j b_m^{-1} \otimes e_{mm} \Big) \Big\|_w \\ & \le & \Big\| \Big( \sum_{j=1}^2 \Big| \Big( \sum_{m \ge 1} a_m^j \otimes e_{1m} \Big)^* \Big|^2 \Big)^\frac12 \Big\|_w \Big\| \Big( \sum_{j=1}^2 \Big| \sum_{m \ge 1} b_m^j b_m^{-1} \otimes e_{mm} \Big|^2 \Big)^\frac12 \Big\|_\infty \\ & \le & \Big\| \Big( \sum_{j=1}^2 \Big| \Big( \sum_{m \ge 1} a_m^j \otimes e_{1m} \Big)^* \Big|^2 \Big)^\frac12 \Big\|_w \ = \ \Big\| \sum_{j=1}^2 \Big| \Big( \sum_{m \ge 1} a_m^j \otimes e_{1m} \Big)^* \Big|^2 \Big\|_{\frac{w}{2}}^{\frac12}.
\end{eqnarray*}
The desired inequality follows them for the fact that $\| \cdot \|_{w/2}$ is a $\frac{w}{2}$-norm for $w < 2$. The expected upper bound for $B$ is easier to obtain since $s > 2$ and we may use the triangle inequality in $L_{s/2}(\M)$ $$B = \Big\| \sum_{j=1}^2 \sum_{m \ge 1} |b_m^j|^2 \Big\|_{\frac{s}{2}}^\frac12 \le \Big( \sum_{j=1}^2 \Big\| \sum_{m \ge 1} b_m^j \otimes e_{m1} \Big\|_s^2 \Big)^{\frac12} \le \Big( \sum_{j=1}^2 \Big\| \sum_{m \ge 1} b_m^j \otimes e_{m1} \Big\|_s^w \Big)^{\frac1w}.$$ This proves ii). The proof of i) is simpler since $L_{w/2}(\M)$ is a Banach space. \fin 

\subsection{$L_pmo$ spaces and duality}

We shall need in what follows to consider the duals of the spaces considered so far. Let $2 \le p' \le \infty$ and $s \ge 2$ so that $w'$ is given by $1/w' = 1/p' + 1/s$. Then we define 
\begin{eqnarray*}
L_{p'w'}^c mo(\M) & = & \Big\{ x = (x_n)_{n \ge 1} \mbox{ $L_{w'}$-martingale} \, : \ \|x\|_{L_{p'w'}^cmo} < \infty \Big\}, \\
L_{p'w'}^{1_c}mo(\M) & = & \Big\{ x = (x_n)_{n \ge 1} \mbox{ $L_{w'}$-martingale} \, : \ \|x\|_{L_{p'w'}^{1_c}mo} < \infty \Big\},  
\end{eqnarray*}
where
\begin{eqnarray*}
\|x\|_{L^c_{p'w'}mo} \hskip-6pt & = & \hskip-8pt \sup_{\begin{subarray}{c} m \ge 1 \\ b_n \in L_s(\M_n) \\ \| \sum_n b_n b_n^* \|_{\frac{s}{2}} \le 1 \end{subarray}} \hskip-3pt \Big\| \Big( \sum_{n \le m} (x_m - x_n) b_n b_n^* (x_m - x_n)^* \Big)^\frac12 \Big\|_{w'}, \\ 
\|x\|_{L^{1_c}_{p'w'}mo} \hskip-6pt & = & \hskip-8pt \sup_{\begin{subarray}{c} \beta_n \in L_s(\M) \\ \| \sum_n \beta_n \beta_n^* \|_{\frac{s}{2}} \le 1 \end{subarray}} \hskip-3pt \Big\| \Big( \sum_{n \ge 1} (x_n - x_{n-1}) \beta_n \beta_n^* (x_n - x_{n-1})^* \Big)^\frac12 \Big\|_{w'}.
\end{eqnarray*}
As usual, we take adjoints to define the row spaces. We should also recall that $L_{p'2}^cmo(\M)$ coincides with the $L_{p'}^cmo(\M)$ spaces introduced in \cite{Per0}. We are now proving that Fefferman's $\mathrm{H}_1-\mathrm{BMO}$ duality theorem extends to these spaces. 

\begin{lemma} \label{pro:fefferman stein duality}
If $1<p<2$ and $w \ge2$, we find $$h^c_{pw}(\M)^*\simeq L^c_{p'w'}mo(\M) \quad \mbox{and} \quad h^{1_c}_{pw}(\M)^* \simeq L^{1_c}_{p'w'}mo(\M).$$
In addition, the analogous duality results also hold for the corresponding row spaces.
\end{lemma}

\vskip-2pt

\dem Again we only consider the column cases. Let us first study the duality for $h_{pw}^c(\M)$. Let $(x,y) \in L_{p'w'}^cmo(\M) \times h^c_{pw}(\M)$. Given any $\delta>0$, we may find a decomposition $y = \sum_{n \ge 1} a_n b_n$ with $\E_n(a_n)=0$, $b_n \in L_s(\M_n)$ and $$\Big\| \sum_{n \ge 1} a_n \otimes e_{1n} \Big\|_w \Big\| \sum_{n \ge 1} b_n \otimes e_{n1} \Big\|_s \le (1+\delta) \|y\|_{h^c_{pw}(\M)}.$$ Then we have
\begin{eqnarray*}
| \langle x, y \rangle| & \le & \sup_{m \ge 1} \Big| \sum_{n \ge 1} \tau \big( x_m^* a_n b_n \big) \Big| \\ & = & \sup_{m \ge 1} \Big| \sum_{n \le m} \tau \big( (x_m-x_n)^* a_n b_n \big) \Big| \\ & \le & \sup_{m \ge 1} \Big\| \Big( \sum_{n \le m} (x_m-x_n)b_n^* b_n (x_m -x_n)^* \Big)^\frac12 \Big\|_{w'} \Big\| \sum_{n \ge 1} a_n \otimes e_{1n} \Big\|_w \\ & \le & \|x\|_{L^c_{p'w'}mo} \Big\| \sum_{n \ge 1} a_n \otimes e_{1n} \Big\|_w \Big\| \sum_{n \ge 1} b^*_n \otimes e_{1n} \Big\|_s \lesssim_\delta \|x\|_{L^c_{p'w'}mo} \|y\|_{h^c_{pw}(\M)}.
\end{eqnarray*}
This proves that the map $$\Phi: L_{p'w'}^cmo(\M) \ni x \mapsto \Phi_x \in h_{pw}^c(\M)^*$$ given by $\Phi_x(y) = \sum_k \tau(d_k(x) d_k(y))$ is bounded. To justify that $\Phi$ is an embedding let us prove that some ball of $h^c_{pw}(\M)$ is norming in ${L^c_{p'w'}mo}$. Namely, if we let $L_p(\M;\ell_2^r)$ denote the space of sequences in $L_p(\M)$ with norm $\|\sum_n (a_n a_n^*)^{1/2} \|_p$ and we write $L_p^{\mathrm{ad}}(\M;\ell_2)$ for the subspace of adapted sequences 
\begin{eqnarray*}
\|x\|_{{L^c_{p'w'}mo}} & = & \sup_{m \ge 1} \Big\{ \big\| \big( (x_m-x_n)b_n \big)_{n \ge 1} \big\|_{L_{w'}(\ell^r_2)} : \|b\|_{L^{ad}_s(\ell^r_2)} \le 1 \Big\} \\
& = & \sup_{m \ge 1} \Big\{ \Big| \sum_{n \ge 1} \tau \Big( \big( (x_m-x_n)b_n \big)^*\eta_n \Big) \Big| : \|b\|_{L^{\mathrm{ad}}_s(\ell^r_2)}, \|\eta\|_{L_w(\ell^r_2)} \le 1 \Big\} \\
& = & \sup_{m \ge 1} \Big\{ \Big| \sum_{n \ge 1} \tau \Big( x_m^* \big( \eta_n-\E_n(\eta_n) \big) b^*_n \Big) \Big| : \|b\|_{L^{\mathrm{ad}}_s(\ell^r_2)}, \|\eta\|_{L_w(\ell^r_2)} \le 1 \Big\} \\ & \le & \sup_{m \ge 1} \Big\{ \big| \tau(x_m^*y) \big| \, : \|y\|_{h^c_{pw}(\M)} \le c_w \Big\} \ = \ c_w \sup_{\|y\|_{h_{pw}^c(\M)} \le 1} |\Phi_x(y)|. 
\end{eqnarray*}
The last inequality follows from the noncommutative Stein inequality \cite{PX1} since $$\Big\| \sum_{n \ge 1} \big( \eta_n - \E_n(\eta_n) \big) b^*_n \Big\|_{h^c_{pw}(\M)} \le \Big\| (\eta_n-\mathcal{E}_n(\eta_n))_n \Big\|_{L_w(\ell^r_2)} \|b\|_{L_s(\ell^r_2)} \le c_w \|\eta\|_{L_w(\ell^r_2)}.$$ This proves that $\Phi$ is indeed an embedding. To prove it is surjective it suffices to show that every continuous functional in $h_{pw}^c(\M)$ is of the form $\Phi_x$ for some $L_{w'}$-martingale $x$ and use the inequality above to justify $x \in L_{p'w'}^cmo(\M)$. The trivial inclusion $h_{pw}^c(\M) \subset L_p(\M)$ shows that every such functional is of the usual form $y \mapsto \tau(z^*y)$ for some $z \in L_{p'}(\M)$. Therefore, $x$ is given by $x_n = \E_n(z)$ which is an $L_{p'}$-martingale and thus and $L_{w'}$-martingale since $w' < p'$.

The duality for $h_{pw}^{1_c}(\M)$ is similar. Let $x \in L_{p'w'}^{1_c}(\M)$ and $y \in h^{1_c}_{pw}(\M)$. By definition, given any $\delta >0$ we may assume that there exists a decomposition of $y = \sum_{n \ge 1} d_n(\alpha_n \beta_n)$ such that $$\Big\| \sum_{n \ge 1} \alpha_n \otimes e_{1n} \Big\|_w \Big\| \sum_{n \ge 1} \beta_n \otimes e_{n1} \Big\|_s \le (1+\delta) \|y\|_{h^{1_c}_{pw}(\M)}.$$ Then we have
\begin{eqnarray*}
| \langle x, y \rangle | \hskip-5pt & = &\hskip-5pt \Big| \sum_{n \ge 1} \tau \big( d_n(x)^* \alpha_n \beta_n \big) \Big| \\ \hskip-5pt & \le &\hskip-5pt \Big\| \Big( \sum_{n \ge 1}d_n(x)\beta^*_n \beta_n d_n(x)^* \Big)^\frac12 \Big\|_{w'} \Big\| \sum_{n \ge 1} \alpha_n \otimes e_{1n} \Big\|_w \\
\hskip-5pt & \le &\hskip-5pt \|x\|_{L^{1_c}_{p'w'}mo} \Big\| \sum_{n \ge 1} \alpha_n \otimes e_{1n} \Big\|_w \Big\| \sum_{n \ge 1} \beta_n \otimes e_{n1} \Big\|_s \lesssim_\delta \|x\|_{L^{1_c}_{p'w'}mo} \|y\|_{h^{1_c}_{pw}(\M)}.
\end{eqnarray*}
For the reverse inequality we note that
\begin{eqnarray*}
\|x\|_{{L^{1_c}_{p'w'}mo}} & = &\sup \Big\{ \big\| (d_n(x) \beta_n)_n\|_{L_{w'}(\ell^r_2)} : \|\beta\|_{L_s(\ell^r_2)} \le 1 \Big\} \\
& = & \sup \Big\{ \Big| \sum_{n \ge 1} \tau \big( (d_n(x) \beta_n)^*\eta_n \big) \Big| : \|\beta\|_{L_s(\ell^r_2)}, \|\eta\|_{L_w(\ell^r_2)} \le 1 \Big\} \\
& = & \sup \Big\{ \Big| \sum_{n \ge 1} \tau \big( d_n(x)^*d_n(\eta_n \beta^*_n) \big) \Big| : \|\beta\|_{L_s(\ell^r_2)}, \|\eta\|_{L_w(\ell^r_2)} \le 1 \Big\} \\
& \le & \sup_{\|y\|_{h^{1_c}_{pw}(\M)} \le 1} |\Phi_x(y)|, \quad \mbox{where} \quad \Phi_x(y) = \summ_k \tau \big(d_k(x) d_k(y) \big).
\end{eqnarray*}
Surjectivity follows again from $h_{pw}^{1_c}(\M) \subset L_p(\M)$, see Lemma \ref{lem:hpw to lp} for the proof. \fin

\subsection{Proof of Theorem Bi}

We now turn to the proof of Theorem Bi, for which we recall the definition of the spaces $\H_{pw}^r(\M)$ and $\H_{pw}^c(\M)$. Given $1 < p < 2$ and $w \ge 2$, we set $$\H_{pw}^r(\M) = h_{pw}^r(\M) + h_{pw}^{1_r}(\M) \quad \mbox{and} \quad \H_{pw}^c(\M) = h_{pw}^c(\M) + h_{pw}^{1_c}(\M).$$ Recall that Theorem Bi for $w=2$ follows from the Burkholder-Davis type inequality in \cite{JuPe14} and Remark \ref{w=2}. We may therefore assume in what follows that $w >2$. Let us start with the inclusion 
$$\H^r_{pw}(\M) + \H^c_{pw}(\M) \subset L_p(\M),$$ which holds for a wider range of $p$'s and $w$'s, as we justify in the following result.

\begin{lemma}\label{lem:hpw to lp}
The continuous inclusion $$\H^r_{pw}(\M) + \H^c_{pw}(\M) \subset L_p(\M)$$ holds for all $1\leq p\leq2$, $w\geq p$ and $s\geq2$ provided $1/p=1/w+1/s$. 
\end{lemma}

\dem It suffices to prove the continuous inclusion of the column spaces, which in turn reduces to prove it for $h^c_{pw}(\M)$ and $h_{pw}^{1_c}(\M)$. The first space embeds trivially in $L_p(\M)$ from H\"older inequality. The second embedding is more involved and we shall divide the proof into three cases:

\noindent 1. \emph{The case $w\ge2$}. According to the definition of the space $h_{pw}^{1_c}(\M)$, we may write $x$ in the form $\sum_n d_n(\alpha_n \beta_n)$ with $\alpha, \beta$ being sequences in $L_w(\M)$ and $L_s(\M)$ respectively. Then we may use the factorization identity $\E_k(\alpha \beta) = u_k(\alpha^*)^* u_k(\beta)$ for a right $\M_k$-module map $u_k: \M \to C(\M_k)$ and the noncommutative dual Doob inequality in $L_{w/2}(\M)$ and $L_{s/2}(\M)$ \cite{Jun02} to conclude that
\begin{eqnarray*}
\|x\|_p & = & \Big\| \sum_{n \ge 1} d_n(\alpha_n \beta_n) \Big\|_p \\ & \le & \Big\| \sum_{n \ge 1} \E_n(\alpha_n \beta_n) \Big\|_p + \Big\| \sum_{n \ge 2} \E_{n-1}(\alpha_n \beta_n) \Big\|_p \\ & \le & \Big\| \sum_{n \ge 1} u_n(\alpha_n^*)^* \otimes e_{1n} \Big\|_w \Big\| \sum_{n \ge 1} u_n(\beta_n) \otimes e_{n1} \Big\|_s \\ & + & \Big\| \sum_{n \ge 2} u_{n-1}(\alpha_n^*)^* \otimes e_{1n} \Big\|_w \Big\| \sum_{n \ge 2} u_{n-1}(\beta_n) \otimes e_{n1} \Big\|_s \\ & \le & 2 \Big\| \Big( \sum_{n \ge 1} |\alpha_n^*|^2 \Big)^\frac12 \Big\|_w \Big\| \Big( \sum_{n \ge 1} |\beta_n|^2 \Big)^\frac12 \Big\|_s \ \le \ 2 \|x\|_{h_{pw}^{1c}(\M)}.
\end{eqnarray*}
Note the last inequality follows by taking infimums over $\alpha$ and $\beta$ as above. 

\noindent 2. \emph{The case $w=p$.} The Burkholder-Gundy and Stein inequalities from \cite{PX1} yield
\begin{eqnarray*}
\|x\|_p & = & \Big\| \sum_{n \ge 1} d_n(\alpha_n \beta_n) \Big\|_p \\ & \le & c_p \Big\| \sum_{n \ge 1} d_n(\alpha_n \beta_n) \otimes e_{1n} \Big\|_p \\ & \le & c_p \Big\| \sum_{n \ge 1} \alpha_n \beta_n \otimes e_{1n} \Big\|_p \ = \ c_p \Big\| \sum_{n \ge 1} \alpha_n \beta_n \beta_n^* \alpha^*_n \Big\|^{\frac12}_{\frac{p}{2}} \\ & \le & c_p \Big\| \sum_{n \ge 1} \alpha_n \otimes e_{1n} \Big\|_p \Big\| \sum_{n \ge 1} \beta_n \otimes e_{n1} \Big\|_\infty \ \le \ c_p \|x\|_{h_{pw}^{1_c}(\M)}.
\end{eqnarray*}
Again, the last inequality follows by taking infimums since $1/s = 1/p - 1/w = 0$.

\noindent 3. \emph{The case $p < w < 2$.} Note that we may assume $p>1$, since we have $w=s=2$ for $p=1$. We proceed by complex interpolation. Let $0 < \theta < 1$ be determined by $1/w=(1-\theta)/p+\theta/2$ and then fix $r = \theta s$. Let $\partial_j$ be the vertical line in $\C$ of complex numbers $z$ with $\mathrm{Im}(z)=j$ for $j=0,1$. Then, we can find two sequences of operator-valued analytic functions $A(z) = (\alpha_n(z))_{n \ge 1}$ and $B(z) = (\beta_n(z))_{n \ge 1}$ satisfying $(A(\theta), B(\theta)) = (\alpha_n,\beta_n)_{n \ge1}$ and 
\begin{eqnarray*}
\max \Big\{ \sup_{z \in \partial_0} \Big\| \sum_{n \ge 1} \alpha_n(z) \otimes e_{1n} \Big\|_p \hskip2pt , \sup_{z \in \partial_1} \Big\| \sum_{n \ge 1} \alpha_n(z) \otimes e_{1n} \Big\|_2 \Big\} \hskip-5pt & \le & \hskip-5pt \Big\| \sum_{n \ge 1} \alpha_n \otimes e_{1n}  \Big\|_w, \\ \max \Big\{ \sup_{z \in \partial_0} \Big\| \sum_{n \ge 1} \beta_n(z) \otimes e_{n1} \Big\|_\infty, \sup_{z \in \partial_1} \Big\| \sum_{n \ge 1} \beta_n(z) \otimes e_{n1} \Big\|_r \Big\} \hskip-5pt & \le & \hskip-5pt \Big\| \sum_{n \ge 1} \beta_n \otimes e_{n1} \Big\|_s \hskip2pt .
\end{eqnarray*}
Note that $r \ge 2$ since $$\frac{\theta}{r} = \frac{1}{s} = \frac{1}{p} - \frac{1}{w} = \frac{\theta}{p} - \frac{\theta}{2} \ \Rightarrow \ r = \frac{2p}{2-p} \ge 2.$$ Then, the three lines lemma and the previous two cases give rise to
\begin{eqnarray*}
\|x\|_p & \le & \sup_{z_j \in \partial_j} \Big\| \sum_{n \ge 1} d_n(\alpha_n(z_0) \beta_n(z_0)) \Big\|^{1-\theta}_p \Big\| \sum_{n \ge 1} d_n(\alpha_n(z_1) \beta_n(z_1)) \Big\|^{\theta}_{p} \\ & \le & c_p \sup_{z \in \partial_0} \Big\| \sum_{n \ge 1} \alpha_n(z) \otimes e_{1n} \Big\|_p^{1-\theta} \Big\| \sum_{n \ge 1} \beta_n(z) \otimes e_{n1} \Big\|_\infty^{1-\theta} \\ & \times & 2^\theta \sup_{z \in \partial_1} \Big\| \sum_{n \ge 1} \alpha_n(z) \otimes e_{1n} \Big\|_2^\theta \Big\| \sum_{n \ge 1} \beta_n(z) \otimes e_{n1} \Big\|_r^\theta \\ & \le & c_p \Big\| \sum_{n \ge 1} \alpha_n \otimes e_{1n} \Big\|_w \Big\| \sum_{n \ge 1} \beta_n \otimes e_{n1} \Big\|_s,
\end{eqnarray*}
which implies the assertion by taking infimums over $\alpha$ and $\beta$ as usual. \fin

The following lemma is a dual version of Lemma \ref{lem:hpw to lp}. It will be used below in an extrapolation argument to obtain the remaining embedding for Theorem Bi though the duality result established in Lemma \ref{pro:fefferman stein duality}.

\begin{lemma}\label{lem:lp to dual of hpw}
The continuous inclusion $$L_{p'}(\M) \subset L^r_{p'w'}mo(\M) \cap L^{1_r}_{p'w'}mo(\M) \cap L^c_{p'w'}mo(\M) \cap L^{1_c}_{p'w'}mo(\M)$$ holds for all $2 < p' < \infty$, $1 < w' \le p'$ and any $s\geq2$ provided $1/w'=1/p'+1/s$. 
\end{lemma}

\dem Arguing as in the proof of Lemma \ref{pro:fefferman stein duality} and using Lemma \ref{lem:hpw to lp}
\begin{eqnarray*}
\|x\|_{{L^c_{p'w'}mo}} & \le & \sup \Big\{ |\langle x, y \rangle | : \|y\|_{h^c_{pw}} \le c_w \Big\} \\
& \le & \sup \Big\{ \|x\|_{p'} \|y\|_p : \|y\|_{h^c_{pw}} \le c_w \Big\} \ \le \ c_{pw} \|x\|_{p'}.
\end{eqnarray*}
Similar estimates hold for $L^{1_c}_{p'w'}mo(\M)$, as well as for their row analogues. \fin

\begin{lemma}\label{lem:dual of hpw to lp}
The continuous inclusion $$L^r_{p'w'}mo(\M) \cap L^{1_r}_{p'w'}mo(\M) \cap L^c_{p'w'}mo(\M) \cap L^{1_c}_{p'w'}mo(\M) \subset L_{p'}(\M)$$ holds for all $2 < p' < \infty$, $1 < w' \le 2$ and any $s\geq2$ provided $1/w'=1/p'+1/s$. 
\end{lemma}

\dem The case $w'=2$ follows from Lemma \ref{pro:fefferman stein duality} and the Davis type inequality proved in \cite{JuPe14}. Let us then fix $1 < w' < 2$ for what follows. We may assume that $x=\mathcal{E}_m(x)$ is a finite martingale and prove the result with constants independent of $m$. Note that for $w'<2$, we can choose $\widetilde{w}', \widetilde{s} > 2$ such that $$\frac{1}{\widetilde{w}'}-\frac{1}{\widetilde{s}}=\frac1p=\frac{1}{{w}'}-\frac{1}{{s}}.$$ This means there exists some $0<\theta<1$ satisfying
$$\frac{1}{2} = \frac{1-\theta}{w'}+\frac{\theta}{\widetilde{w}'} \ \Rightarrow \ \frac1q := \frac12-\frac 1p = \Big( \frac{1-\theta}{w'}+\frac{\theta}{\widetilde{w}'} \Big) - \Big( \frac{1}{{w}'}-\frac{1}{{s}} \Big) = \frac{1-\theta}{s}+\frac{\theta}{\tilde{s}}.$$
Let $(b_n)_{n \ge 1}$ be an adapted sequence in $L_q(\M)$ satisfying $\|\sum_n b_n b_n^*\|_q \le 1$. Then we can find a sequence of vector-valued analytic functions $B(z) = (b_n(z))_{n \ge 1}$ with $$B(\theta) = (b_n)_{n \ge 1} \quad \mbox{and} \quad \max \Big\{ \sup_{z \in \partial_0} \Big\| \sum_{n \ge 1} b_n(z) \otimes e_{1n} \Big\|_s, \sup_{z \in \partial_1} \Big\| \sum_{n \ge 1} b_n(z) \otimes e_{1n} \Big\|_{\widetilde{s}} \Big\} \le 1.$$ Thus, we deduce from Lemma \ref{lem:lp to dual of hpw} that
\begin{eqnarray*}
\lefteqn{\hskip-15pt \Big\| \Big( \sum_{n \le m} (x_m - x_n) b_n b_n^* (x_m - x_n)^* \Big)^\frac12 \Big\|_2} \\ & \le & \sup_{z \in \partial_0} \Big\| \Big( \sum_{n \ge 1} (x_m - x_n) b_n(z) b_n(z)^* (x_m - x_n) \Big)^\frac12 \Big\|_{w'}^{1-\theta} \\ & \times & \sup_{z \in \partial_1} \Big\| \Big( \sum_{n \ge 1} (x_m - x_n) b_n(z) b_n(z)^* (x_m - x_n)^* \Big)^\frac12 \Big\|_{\widetilde{w}'}^{\theta} \\ [5pt] & \le & \|x\|_{L^c_{p'w'}mo}^{1-\theta} \|x\|_{L^c_{p'\widetilde{w}'}mo}^{\theta} \ \le \ c_{p\widetilde{w}'} \|x\|_{L^c_{p'w'}mo}^{1-\theta} \|x\|_p^{\theta}.
\end{eqnarray*}
According to the definition of $L_{p'w'}^cmo(\M)$, we immediately conclude that it embeds in $L_{p'}(\M)$. The exact same argument applies for $L_{p'w'}^{1_c}mo(\M)$ and row spaces. \fin

\begin{remark}\label{rem:JohnNirenberg}
\emph{Lemmas \ref{lem:lp to dual of hpw} and \ref{lem:dual of hpw to lp} yield a John-Nirenberg type result for $p>2$.}
\end{remark}

Applying the Duality Lemma \ref{pro:fefferman stein duality} to our embedding in Lemma \ref{lem:dual of hpw to lp}, we obtain the converse embedding of Lemma \ref{lem:hpw to lp}. Altogether, this proves the first assertion in Theorem Bi. The second assertion follows from the following result. 

\begin{lemma}\label{lem:hcpw to Hcp}
If $1<p<2$ and $w \ge 2$, we find $$\H_{pw}^c(\M) \subset \H_p^c(\M)$$ up to a constant $c_{pw}$. The same continuous inclusions hold in the row case. 
\end{lemma}

\dem
We shall prove that
\begin{eqnarray*}
\|x\|_{\H^c_p(\M)} & \le & c_{pw} \|x\|_{h^c_{pw}(\M)}, \\
\|x\|_{\H^c_p(\M)} & \le & c_{pw} \|x\|_{h^{1_c}_{pw}(\M)}.
\end{eqnarray*}
For the first inequality, assume that $x = \sum_n a_n b_n$ with $$\Big\| \sum_{n \ge 1} a_n \otimes e_{1n} \Big\|_w \Big\| \sum_{n \ge 1} b_n \otimes e_{n1} \Big\|_s \le (1+\delta) \|x\|_{h^c_{pw}(\M)},$$ where $\E_{n}(a_n)=0$ and $b_n \in L_s(\M_n)$. By Davis decomposition \cite{Per}, we have
\begin{eqnarray*}
\|x\|_{\H^c_p(\M)} & \lesssim & \|x\|_{h^c_p(\M)} \\ [3pt] & = & \Big\| \Big( \sum_{n \ge 1} \E_{n-1}|d_n(x)|^2 \Big)^{\frac12} \Big\|_p \\ & = & \Big\| \Big( \sum_{n \ge 1} \E_{n-1} \Big|\sum_{k<n} d_n(a_kb_k) \Big|^2 \Big)^{\frac12} \Big\|_p \\ & = & \Big\| \Big( \sum_{n \ge 1} \E_{n-1} \Big| \sum_{k<n} d_n(a_k)b_k \Big|^2 \Big)^{\frac12} \Big\|_p.
\end{eqnarray*}
By the right $\M_{n-1}$-modular maps $u_{n-1}: \M \to C(\M_{n-1})$ from \cite[Proposition 2.8]{Jun02}
\begin{eqnarray*}
\mathcal{E}_{n-1} \Big| \sum_{k<n} d_n(a_k)b_k \Big|^2 & = & \Big| u_{n-1} \Big( \sum_{k<n} d_n(a_k)b_k \Big) \Big|^2 \\
& = & \Big| \sum_{k<n} u_{n-1}(d_n(a_k))b_k \Big|^2 \ = \ \big| \widehat{u}_{n-1}(\widehat{d}_n(A))B \big|^2,
\end{eqnarray*}
where $A = \sum_k a_k \otimes e_{1k}, B = \sum_k b_k\otimes e_{k1}$ and $\widehat{m}$ is used for the matrix amplification $m \otimes id_{\mathcal{B}(\ell_2)}$ of the map $m$. Then, use H\"older inequality and Burkholder-Gundy inequality in the case $w>2$
\begin{eqnarray*}
\|x\|_{\H^c_p(\M)} \hskip-5pt & \lesssim & \hskip-5pt \Big\| B^* \sum_{n \ge 1} \big| \widehat{u}_{n-1}(\widehat{d}_n(A)) \big|^2 B \Big\|^{\frac12}_{\frac{p}{2}} \\
\hskip-5pt & \le & \hskip-5pt \Big\| \Big( \sum_{n \ge 1} \widehat{\E}_{n-1} \big| \widehat{d}_n(A) \big|^2 \Big)^{\frac12} \Big\|_w \|B\|_s \le c_w \|A\|_w \|B\|_s \lesssim_\delta c_w \|x\|_{h^c_{pw}(\M)}.
\end{eqnarray*}
For the second inequality, assume that $x = \sum_n d_n(\alpha_n \beta_n)$ with $$\Big\| \sum_{n \ge 1} \alpha_n \otimes e_{1n} \Big\|_w \Big\| \sum_{n \ge 1} \beta_n \otimes e_{n1} \Big\|_s \le (1+\delta) \|x\|_{h^{1_c}_{pw}(\M)}.$$ Now recall that $L_q(\M;\ell_q) \subset L_q(\M;\ell^c_2)$ for $q \le 2$ and the reverse embedding holds for $q \ge 2$. Indeed, the cases $q=2$ and $q=\infty$ are clear. Then one can proceed by interpolation and duality. Thus, noting that $w,s \ge 2$ and $p \le 2$
\begin{eqnarray*}
\hskip17pt \|x\|_{\H^c_p(\M)} & = & \big\| (d_n(\alpha_n \beta_n)_{n \ge 1} \big\|_{L_p(\M; \ell^c_2)} \\ & \le & \big\| (d_n(\alpha_n \beta_n)_{n \ge 1} \big\|_{L_p(\M; \ell_p)} \\ & \le & 2 \big\| (\alpha_n)_{n \ge 1} \big\|_{L_w(\M;\ell_{w})} \big\| (\beta_n)_{n \ge 1} \big\|_{L_s(\M;\ell_s)} \\ & \le &  2 \big\| (\alpha_n)_{n \ge 1} \big\|_{L_w(\M;\ell^r_{2})} \big\|(\beta_n)_{n \ge 1} \big\|_{L_s(\M;\ell^c_2)} \ \lesssim_\delta \ 2 \|x\|_{h_{pw}^{1_c}(\M)}. \hskip18pt \square
\end{eqnarray*}

\subsection{Proof of Theorem Bii}

Our aim now is to complete the proof of Theorem B. The last assertion follows trivially from the Davis type decomposition in Theorem Bi and the boundedness of $(\E_n)_{n \ge 1}: \H_{pw}^\dag(\M) \to L_p(\M;\ell_\infty^\dag)$ for $1 < p < 2$, $w>2$ and $\dag = r,c$. This latter estimate will be our goal. As usual, we only justify the column case. It suffices to show that
\begin{eqnarray*}
\big\| (\E_n(x))_{n \ge 1} \big\|_{L_p(\M; \ell^c_\infty)} & \le & c_{pw} \|x\|_{h^c_{pw}(\M)}, \\ 
\big\| (\E_n(x))_{n \ge 1} \big\|_{L_p(\M; \ell^c_\infty)} & \le & c_{pw} \|x\|_{h^{1_c}_{pw}(\M)}.
\end{eqnarray*}
Given $x \in h^c_{pw}(\M)$, there exists a decomposition $x = \sum_{n} a_n b_n$ satisfying $$\Big\| \sum_{n \ge 1} a_n \otimes e_{1n} \Big\|_w \Big\| \sum_{n \ge 1} b_n \otimes e_{n1} \Big\|_s \le (1+\delta) \|x\|_{h^c_{pw}(\M)}$$ where $\E_n(a_n)=0$ and $b_n \in L_s(\M_n)$. This gives rise to $\E_n(x) = \widehat{\E}_n(A)B$ where $A = \sum_k a_k\otimes e_{1k}$ and $B = \sum_k b_k \otimes e_{k1}$ as in the proof of Theorem Ai. Using polar decomposition $A = v_A |A|$ and the modular map $u_n: \M \to C(\M_n)$, we can rewrite $$\widehat{\E}_n(A) = \widehat{\E}_n(v_A|A|) = \widehat{u}_n(v_A^*)^* \widehat{u}_n(|A|).$$ Note that $\widehat{u}_n(v_A^*)^*$ is a contractive row. On the other hand, since $w>2$ we may use noncommutative Doob inequality for $|A|^2 \in L_{w/2}(\M \bar\otimes \mathcal{B}(\ell_2))$. Doing so we deduce there exists $\beta \in L_w(\M \bar\otimes \mathcal{B}(\ell_2))$ and contractions $\gamma_n \in \M \bar\otimes \mathcal{B}(\ell_2)$ satisfying $$\widehat{\E}_n(|A|^2) = \beta^* \gamma_n^* \gamma_n \beta \quad \mbox{and} \quad \|\beta\|_w \le \|A\|_w.$$ This implies $\widehat{u}_n(|A|) = \psi_n \gamma_n \beta$ for some contraction $\psi_n \in \M \bar\otimes \mathcal{B}(\ell_2)$. As we did in the proof of Theorem Ai, we now exploit that $\widehat{u}_n(v_A^*)^*$ is a row and $B$ is a column to find a factorization of $\E_n(x)$ with operators affiliated to $\M$. Namely $$\E_n(x) = \widehat{\E}_n(A)B = \underbrace{\widehat{u}_n(v_A^*)^* \psi_n \gamma_n \beta B (B^* \beta^* \beta B)^{-\frac12}}_{w_n} \underbrace{(B^* \beta^* \beta B)^{\frac12}}_{b}.$$ Since $w_n \in \M$ is a contraction, we just need to observe that 
\begin{eqnarray*}
\|b\|_{L_p(\M)} & = & \|\beta B\|_{L_p(\M \bar\otimes \mathcal{B}(\ell_2))} \\ [3pt] & \le & \|\beta\|_{L_w(\M \bar\otimes \mathcal{B}(\ell_2))} \|B\|_{L_s(\M \bar\otimes \mathcal{B}(\ell_2))} \\ & \le & \Big\| \sum_{n \ge 1} a_n \otimes e_{1n} \Big\|_w \Big\| \sum_{n \ge 1} b_n \otimes e_{n1} \Big\|_s \lesssim_\delta \|x\|_{h^c_{pw}(\M)}.
\end{eqnarray*} 
Let us finally prove the inequality for $x = \sum_n d_n (\alpha_n \beta_n) \in h_{pw}^{1_c}(\M)$ with $$\Big\| \sum_{n \ge 1} \alpha_n \otimes e_{1n} \Big\|_w \Big\| \sum_{n \ge 1} \beta_n \otimes e_{n1} \Big\|_s \le (1+\delta) \|x\|_{h^{1_c}_{pw}(\M)}.$$ Then we can write $$\E_n(x) = \sum^n_{k=1} d_k(\alpha_k \beta_k) = \sum^n_{k=1} \E_k(\alpha_k \beta_k) - \sum^n_{k=1} \E_{k-1}(\alpha_k \beta_k) = X_n - Y_n.$$ Using the quasi-triangle inequality, we are reduced to deal with $(X_n)_{n \ge 1}$ and $(Y_n)_{n \ge 1}$. The two cases being similar, we only estimate $(X_n)_{n \ge 1}$. Then, using the modular map $u_k: \M \to C(\M_k)$ we may write $$X_n = \sum^n_{k=1}u_k(\alpha^*_k)^*u_k(\beta_k) = \Big( \sum^n_{k=1} u_k(\alpha^*_k)^* \otimes e_{1k} \Big) \Big( \sum_{k \ge 1} u_k(\beta_k) \otimes e_{k1} \Big) = \mathbb{E}_n(A)B,$$
where $\mathbb{E}_n = id_{\M} \otimes id_{\mathcal{B}(\ell_2)} \otimes \mathsf{E}_n$ with $\mathsf{E}_n$ being the conditional expectation on $\mathcal{B}(\ell_2)$ which we introduced to prove \eqref{BasicIneqs2}. Here $A = \sum_{k \ge 1} u_k(\alpha_k^*)^* \otimes e_{1k}$. Arguing as above, it all reduces to show that $\|A\|_w \|B\|_s \lesssim_\delta \|x\|_{h_{pw}^{1_c}(\M)}$. However, since $w \ge 2$ we may use the dual form of Doob inequality in $L_{w/2}(\M)$ $$\|A\|_w = \Big\| \sum_{k \ge 1} \E_k |\alpha_k^*|^2 \Big\|^{\frac12}_{\frac{w}{2}} \le c_{\frac{w}{2}} \Big\| \sum_{k \ge 1} |\alpha_k^*|^2 \Big\|^{\frac12}_{\frac{w}{2}} = c_{\frac{w}{2}} \Big\| \sum_{k \ge 1} \alpha_k \otimes e_{1k} \Big\|_w$$ and similarly $\|B\|_s \le c_{s/2} \| \sum_k \beta_k \otimes e_{k1} \|_s$ since $s \ge 2$. The proof is complete. \fin

\begin{remark}
\emph{Since $\H_{pw}^\dag(\M) \subset \H_{p}^\dag(\M)$ for $\dag = r,c$, it turns out from Theorems A and B that other asymmetric Doob maximal inequalities hold for these Hardy spaces. Namely, given $1 < p < 2$ and $w \ge 2$, we have}
\begin{eqnarray*}
(\E_n)_{n \ge 1}: \H_{pw}^r(\M) \to L_p(\M;\ell_\infty^\theta) & \mbox{\emph{for}} & 0 \le \theta < p/2, \\
(\E_n)_{n \ge 1}: \H_{pw}^c(\M) \to L_p(\M;\ell_\infty^\theta) & \mbox{\emph{for}} & 1 - p/2 < \theta \le 1. 
\end{eqnarray*}
\end{remark}

\vskip3pt

\noindent \textbf{Acknowledgement.} Junge is partially supported by the NSF DMS-1201886 and NSF DMS-1501103. Parcet is partially supported by ERC StG-256997-CZOSQP. All the authors are also supported in part by the ICMAT Severo Ochoa Grant SEV-2011-0087 (Spain).

\bibliographystyle{amsplain}

\enlargethispage{2cm}

\vskip20pt

\hfill \noindent \textbf{Guixiang Hong} \\
\null \hfill School of Mathematics and Statistics \\ 
\null \hfill Wuhan University \\ 
\null \hfill Wuhan 430072. China \\
\null \hfill Instituto de Ciencias Matem{\'a}ticas \\ \null \hfill
CSIC-UAM-UC3M-UCM \\ \null \hfill Consejo Superior de
Investigaciones Cient{\'\i}ficas \\ \null \hfill C/ Nicol\'as Cabrera 13-15.
28049, Madrid. Spain \\ \null \hfill\texttt{guixiang.hong@icmat.es}

\vskip2pt

\hfill \noindent \textbf{Marius Junge} \\
\null \hfill Department of Mathematics
\\ \null \hfill University of Illinois at Urbana-Champaign \\
\null \hfill 1409 W. Green St. Urbana, IL 61891. USA \\
\null \hfill\texttt{junge@math.uiuc.edu}

\vskip2pt

\hfill \noindent \textbf{Javier Parcet} \\
\null \hfill Instituto de Ciencias Matem{\'a}ticas \\ \null \hfill
CSIC-UAM-UC3M-UCM \\ \null \hfill Consejo Superior de
Investigaciones Cient{\'\i}ficas \\ \null \hfill C/ Nicol\'as Cabrera 13-15.
28049, Madrid. Spain \\ \null \hfill\texttt{javier.parcet@icmat.es}
\end{document}